\documentclass[a4paper,12pt]{article}
\usepackage{amssymb,amsmath,amsfonts,amsthm}

\numberwithin{equation}{section}
\usepackage[left=.75 in, right=.75 in,top=.75 in, bottom=.75 in]{geometry}

\providecommand{\abs}[1]{\left\vert#1\right\vert}
\providecommand{\norm}[1]{\left\Vert#1\right\Vert}
\providecommand{\pnorm}[2]{\left\Vert#1\right\Vert_{L^{#2}}}
\providecommand{\pnormspace}[3]{\left\Vert#1\right\Vert_{L^{#2}(#3)}}

\providecommand{\Rn}[1]{\mathbb{R}^{#1}}

\def\wstar{\overset{*}{\rightharpoonup}}
\def\nab{\nabla}
\def\nb{\nabla^\bot}
\def\dt{\partial_t}
\def\dph{\partial_{\Phi}}
\def\hal{\frac{1}{2}}
\def\lep{\lambda_\varepsilon}
\def\ep{\varepsilon}

\def\ale{\abs{\log \ep}}

\def\rest{\hskip 1pt{\hbox to 10.8pt{\hfill\vrule height 7pt width 0.4pt depth 0pt\hbox{\vrule height 0.4pt
width 7.6pt depth 0pt}\hfill}}}

\def\evalu{\hskip 1pt{\hbox to 2pt{\hfill \vrule height -6pt width 0.4pt depth0pt}}}

\DeclareMathOperator{\curl}{curl}
\DeclareMathOperator{\diverge}{div}
\DeclareMathOperator{\supp}{supp}
\DeclareMathOperator{\dist}{dist}

\newtheorem{lem}{Lemma}[section]

\newtheorem{prop}[lem]{Proposition}
\newtheorem{thm}[lem]{Theorem}
\newtheorem{remark}[lem]{Remark}
\newtheorem{thm_intro}{Theorem}

\title{Ginzburg-Landau vortex dynamics driven by an applied boundary current}
\author{Ian Tice\footnote{Supported by an NSF
Postdoctoral Research Fellowship}\\
{\small Brown University, Division of Applied Mathematics}\\
{\small 182 George St., Providence, RI 02912}\\
{\small\tt tice@dam.brown.edu} }

\begin{document}

\maketitle

\begin{abstract}
In this paper we study the time-dependent Ginzburg-Landau equations on a smooth, bounded domain $\Omega \subset \Rn{2}$, subject to an electrical current applied on the boundary.  The dynamics with an applied current are non-dissipative, but via the identification of a special structure in an interaction energy, we are able to derive a precise upper bound for the energy growth.  We then turn to the study of the dynamics of the vortices of the solutions in the limit $\ep \rightarrow 0$.  We first consider the original time scale, in which the vortices do not move and the solutions undergo a ``phase relaxation.''  Then we study an accelerated time scale in which the vortices move according to a derived dynamical law.  In the dynamical law, we identify a novel Lorentz force term induced by the applied boundary current.
\end{abstract}

%%%%%%%%%%%%%%%%%%%%%%%%%%%%%%%%%%%%%%%%%%%%%%%%%%%%%%%%%%%%%%%%%%%%%%%%%%%%%
\section{Introduction}
%%%%%%%%%%%%%%%%%%%%%%%%%%%%%%%%%%%%%%%%%%%%%%%%%%%%%%%%%%%%%%%%%%%%%%%%%%%%

%%%%%%%%%%%%%%%%%%%%%%%%%%%%%%%%%%%%%%%%%%%%%%%%%%%%%%%%%%%%%%%%%%%%%%%%%%%%%
\subsection{Formulation of the equations and boundary conditions}
%%%%%%%%%%%%%%%%%%%%%%%%%%%%%%%%%%%%%%%%%%%%%%%%%%%%%%%%%%%%%%%%%%%%%%%%%%%%%

The Ginzburg-Landau free energy functional on a domain $\Omega \subset \Rn{2}$ is defined for a function $u:\Omega \rightarrow \mathbb{C}$, a vector field $A:\Omega \rightarrow \Rn{2}$, and a real parameter $\ep>0$ by 
\begin{equation}\label{gl_energy}
 F_\ep(u,A) = \hal \int_\Omega \abs{\nab_A u}^2 + \abs{\curl{A}}^2 +\frac{(1-\abs{u}^2)^2}{2\ep^2},
\end{equation}
where $\nab_A u:=\nab u - iAu$ is the covariant gradient of $u$.  Ginzburg and Landau introduced their eponymous functional in 1950  \cite{gl} as a  free energy in a phenomenological model of superconductivity.  In this setting $\Omega$ is thought of as a two-dimensional cross section of a sample of superconducting material; we shall assume $\Omega$ is smooth and bounded.  The function $u$ is known as the ``order parameter,'' and models the relative density and phase of superconducting electrons, with $\abs{u} \approx 1$ indicating the  superconducting state and $\abs{u} \approx 0$ indicating the normal state.  The vector field $A$ is the magnetic vector potential and $h := \curl{A}$ is the induced magnetic field strength in $\Omega$.  The vector $h e_3$, which is orthogonal to $\Omega$, is the induced magnetic field.  In the model,  $\ep>0$ is a parameter depending on the material comprising the superconductor.  We are interested in the regime $\ep \ll 1$, which corresponds to so-called ``extreme type-II'' superconductors.

The equations modeling the dynamics of superconductors, derived by  Gor'kov and Eliashberg in 1968 \cite{gor_el}, are the covariant heat flow for the Ginzburg-Landau functional:
\begin{equation}\label{gork_el}
 \begin{cases}
  \dt u + i\Phi u =  \Delta_A u + \frac{u}{\varepsilon^2}(1-\abs{u}^2) &\text{in } \Omega \times \Rn{+} \\
  \sigma(\dt A + \nab \Phi) = \nab^\bot h + (iu,\nab_A u) & \text{in } \Omega \times \Rn{+}. \\
 \end{cases}
\end{equation}
In these equations we have written $\Delta_A u := (\diverge - iA\cdot)\nab_A u$ for the covariant Laplacian and $\nab^\bot h:=(-\partial_2 h, \partial_1 h)$ for the perpendicular gradient of the induced magnetic field.  The function $\Phi:\Omega \rightarrow \Rn{}$ is the electric potential, and  $E:=-(\dt A + \nab \Phi)$ is the induced electric field.  The constant $\sigma >0$ is the conductivity of the superconducting material, and by Ohm's law the quantity $\sigma E$ is the current of normal (i.e. non-superconducting) electrons in $\Omega$.  The vector field $(iu,\nab_A u)$ is known as the supercurrent, and represents the current of electrons in the superconducting state.  Here we have employed the notation $(a,b):=  \Re(a)\Re(b) + \Im(a) \Im(b)$ for $a,b\in\mathbb{C}$ to mean the inner-product with $\mathbb{C}$ identified with $\Rn{2}$, and $(a,X)$ for $X\in\mathbb{C}^2$ to mean the vector in $\Rn{2}$ with components $(a,X_1)$ and $(a,X_2)$.  

The equations \eqref{gork_el} give rise to the Maxwell equations in $\Omega$.  The second equation in \eqref{gork_el} is Amp\`{e}re's law with the total current given as the sum of the normal and supercurrents.  Faraday's law of induction is seen in the relation $\curl{E} = -\curl(\dt A + \nab \Phi) = -\dt h$.  Taking the inner-product $(iu,\cdot)$ with the first equation in \eqref{gork_el} and taking the divergence of the second, we find $\diverge{E} = -(iu,\dph u)/\sigma$, where we have written $\dph u = \dt u + i\Phi u$.  Then $-(iu,\dph u)/\sigma$ is the charge in $\Omega$, and Gauss's law holds.  The remaining Maxwell equation, Gauss's law of magnetism, follows since $\diverge(h e_3) = \partial_3 h =0.$

The evolution of the superconductor is coupled to the electromagnetic fields on the exterior of the domain,  $\Omega^c = \Rn{2}\backslash \Omega$, by assuming that the electromagnetic fields satisfy Maxwell's equations everywhere \cite{ gor_el, chapman_1,  du_gray, du_3}.  In the absence of surface charges, this gives rise (cf. \cite{griffiths}) to boundary jump conditions coupling the external electric and magnetic fields, $E_{ex}$ and $H_{ex}$, to the fields in $\Omega$.  These are
\begin{equation}\label{phsx_conds}
E\cdot \nu =  E_{ex} \cdot \nu  \text{ and } h = H_{ex},
\end{equation}
where $\nu$ is the outward unit normal on $\partial \Omega$.  We will assume that the exterior fields satisfy the static Maxwell equations, which in particular requires
\begin{equation}\label{maxwell}
 \nb H_{ex} =  - I_{ex},
\end{equation}
where $I_{ex}:\Omega^c \rightarrow \Rn{2}$ is the smooth, time-independent current of normal electrons on the exterior.  If $\Omega^c$ is a conductor, we will assume it has the same conductivity as $\Omega$ so that $I_{ex} = \sigma E_{ex}$.  If $\Omega^c$ is not a conductor, then 
$I_{ex} =0$ and we take $H_{ex}>0$;  $E_{ex}$ need not vanish, but it is frequently assumed to when $I_{ex}$ does.  %The assumption that the external fields are static is limiting; a more realistic model would employ the dynamic Maxwell equations on $\Omega^c$.  Nevertheless, the case of static external fields is a reasonable starting point for the mathematical analysis, and we defer an analysis of dynamic coupling to future work.

It is mathematically convenient to recast the boundary conditions for $E\cdot \nu$ in \eqref{phsx_conds} in terms of a condition for $\nab_A u \cdot \nu$.  Since $\nab^\bot h \cdot \nu  = \nab h \cdot \tau =- \partial_\tau h$ with $\tau$ the unit tangent, we may plug  $h=H_{ex}$ and \eqref{maxwell} into the second equation in \eqref{gork_el} to see 
\begin{equation}\label{compatibility}
(iu,\nab_A u \cdot \nu)  = I_{ex} \cdot \nu - \sigma E \cdot \nu.
\end{equation}
When $\sigma E_{ex} = I_{ex}$ (with $I_{ex}$ possibly $0$) we may then take 
$\nab_A u \cdot \nu=0$, which implies the appropriate boundary condition, $E_{ex}\cdot \nu = E \cdot \nu$.

We will actually study a generalization of these boundary conditions, which for clarity we record now along with the evolution equations:
\begin{equation}\label{tdgl}
  \begin{cases}
   \dt u + i \Phi u = \Delta_A u + \frac{u}{\varepsilon^2}(1-\abs{u}^2) &\text{in } \Omega \times \Rn{+} \\
   \dt A + \nab \Phi = \nab^\bot h + (iu,\nab_A u) & \text{in } \Omega \times \Rn{+} \\
  \nabla_A u \cdot \nu = i u J_{ex} \cdot \nu &\text{on } \partial \Omega \times \Rn{+} \\
  h = H_{ex} & \text{on } \partial \Omega \times \Rn{+} \\
  (u,A,\Phi)\rvert_{t=0} = (u_0,A_0,\Phi_0).
 \end{cases}
\end{equation}
Here $H_{ex}$ is again given by \eqref{maxwell}, but now we take $J_{ex}:\partial \Omega \rightarrow \Rn{2}$ to be any smooth vector field.  To reduce notational clutter, we have assumed  that $\sigma =1$, but all of our results may be modified to handle any fixed $\sigma>0$.  The introduction of $J_{ex}$ is justified in three ways.  First, from a mathematical point of view, the case $J_{ex}\neq 0$ is just a generalization of the case $J_{ex}=0$.  The methods we develop to handle $-\nb H_{ex}=I_{ex}\neq 0$ also handle $J_{ex}\neq 0$, which justifies referring to $J_{ex}$ as a sort of current.  Second, the actual physical jump condition across $\partial \Omega$ is \cite{griffiths}
\begin{equation}
 E\cdot \nu = E_{ex} \cdot \nu + q,
\end{equation}
where $q$ is the surface charge accumulated on $\partial \Omega$.  Plugging the generalized condition $\nab_A u \cdot \nu = i u J_{ex} \cdot \nu$ into \eqref{compatibility} when $I_{ex}=E_{ex}$ yields
\begin{equation}
 E\cdot \nu = E_{ex} \cdot \nu - \abs{u}^2 J_{ex} \cdot \nu,
\end{equation}
which shows that $J_{ex}\cdot \nu$ behaves as a sort of surface charge.  Similarly, when $I_{ex}=0$ but $E_{ex} \neq 0$, which corresponds to an external voltage, \eqref{compatibility} gives
\begin{equation}
 E\cdot \nu = -\abs{u}^2 J_{ex} \cdot \nu
\end{equation}
so that $J_{ex}\cdot \nu$ can also behave as an external voltage would.  This suggests that the generalized boundary condition can be used as an approximate model of surface charge or external voltage.  Third, the currents $J_{ex}$ and $I_{ex}$ are independent, which gives a mechanism for inducing different scales of current forcing.  The reader solely interested in the standard  choice of  boundary conditions $h=H_{ex}$, $\nab_A u \cdot \nu=0$ may simply take $J_{ex}=0$ in all of our analysis.

The novelty in our analysis is in the presence of the applied currents $I_{ex}$ and $J_{ex}$.  Numerous authors \cite{bcps,rub_stern,lin_1,js_2,bos_col,bos_conv,bos,serf_1,serf_2} have studied the non-magnetic analog of the equations \eqref{tdgl} for which an applied current is impossible.  In the magnetic case \eqref{tdgl} has been studied rigorously with  $I_{ex}=J_{ex}=0$ and $H_{ex}$ a constant in \cite{spirn,ss_gamma}.  Several  numerical and formal asymptotic results are available \cite{du_2, du_3, du_gray} when $I_{ex} \neq 0$ and $J_{ex}=0$, a stability analysis of the normal state ($u=0$) with applied current was performed in \cite{almog}, and a 1-D model of a superconducting wire with current was studied in \cite{rub_stern_zum}, but we are aware of no rigorous results in the 2-D magnetic model with applied current or with a surface charge.

%%%%%%%%%%%%%%%%%%%%%%%%%%%%%%%%%%%%%%%%%%%%%%%%%%%%%%%%%%%%%%%%%%%%%%%%%%%%%
\subsection{Definitions and terminology}
%%%%%%%%%%%%%%%%%%%%%%%%%%%%%%%%%%%%%%%%%%%%%%%%%%%%%%%%%%%%%%%%%%%%%%%%%%%%%

We will now record several definitions and bits of terminology that will be used throughout the paper.  For a more thorough exposition of these quantities and of the magnetic Ginzburg-Landau model in general, we refer to the book \cite{ss_book} and the references therein.

The objects of interest in the study of \eqref{tdgl} are the zeroes of the complex-valued function $u$; these are known as vortices.  Each vortex carries an integer topological charge called its degree, which is defined as the winding number of the map  $u/\abs{u}$ on any simple closed curve around the zero.  The energy density (the integrand of $F_\ep$) concentrates around the vortices, with the $i^{th}$ vortex contributing an amount of energy of the order $\pi \abs{d_i} \abs{\log\ep}$ to $F_\ep$, where $d_i$ is the degree of the vortex.  
  
The energy $F_\ep$ possesses a gauge invariance under the pointwise action of the group $\mathbb{U}(1)$: $F_\ep(u,A) = F_\ep(ue^{i\xi},A+\nab\xi)$ for any sufficiently smooth $\xi:\Omega \rightarrow \Rn{}$.  This gauge invariance carries over to solutions of the equations \eqref{tdgl} as well in the sense that if $(u,A,\Phi)$ are solutions, then so are $(ue^{i\xi},A+\nab\xi,\Phi-\dt \xi)$ for  $\xi:\Omega\times\Rn{+} \rightarrow \Rn{}$.  We refer to the change  $u \mapsto ue^{i\xi}$, $A \mapsto A + \nab \xi$, $\Phi \mapsto \Phi - \dt \xi$ as a gauge change.  For solutions to \eqref{tdgl} to be unique, we must eliminate the gauge invariance by ``fixing a gauge.''  In this paper, we work exclusively in what we call the $\Phi=f$ gauge (see Lemma \ref{f_gauge}).  

For any $\ep>0$ and any choice of $J_{ex}$, $H_{ex}$ smooth, the system \eqref{tdgl} is well-posed for all time in any fixed gauge, and the solutions are smooth.  This may be established through a modification of standard results  \cite{chen,du_1}.  For the sake of completeness,  in Appendix \ref{well_posed_section} we make a few brief remarks on the necessary changes and the a priori estimates available for  solutions.  

%Because of the ability to change gauges, the only physically meaningful quantities must be gauge invariant.  The energy $F_\ep$ is gauge invariant, and the modulus $\abs{u}^2$ is obviously gauge invariant as well.  So are the induced magnetic and electric fields, $h=\curl{A}$ and $E= -(\dt A + \nab \Phi)$.  The supercurrent $(iu,\nab_A u)$ and the so-called charge, $(iu,\dph u)$, are as well (here we have written $\dph u := \dt u +i\Phi u$ for the covariant time derivative).  Note that our imposed boundary condition $\nab_A u \cdot \nu = iu J_{ex}\cdot \nu$ is also gauge invariant, which makes it possible for it to have a physical meaning.  

The vortices of a configuration $(u,A)$ are best described through the ``vorticity," $\mu(u,A)$, a gauge-invariant version of the Jacobian determinant of $u$:
\begin{equation}
\mu(u,A) := \curl (iu,\nab_A u)+ \curl A.
\end{equation} 
It is now well-known in the literature (cf. \cite{js} for first results but without magnetic field, and \cite{ss_book, tice} for results with magnetic field) that $\mu(u,A) \approx 2\pi \sum d_i \delta_{a_i},$ where $a_i\in \Omega$, $d_i\in\mathbb{Z}$ are the location and degree of the $i^{th}$ vortex.  Here $\approx$ means close in various norms: $(C^{0,\alpha}(\Omega))^*$, for instance.

As in the non-magnetic case \cite{bbh}, two mechanisms contribute to the energy $F_\ep$.  In the case we are interested in, with $n$ vortices of degree $d_i = \pm 1$, the result is that (roughly speaking)
\begin{equation}\label{en_split_intro}
 F_\ep(u,A) =  n ( \pi  \ale + \gamma)  +  W_{d}(a)  +  O(1).
\end{equation}
The first term on the right is the self-energy of the vortices, which is itself composed of a divergent term $\pi n\ale$ and a finite term $n \gamma$, where $\gamma$ is a known constant related to the structure of a vortex.  The second term is the inter-vortex interaction energy, or magnetic renormalized energy.  For a collection of points $a = (a_1,\dotsc,a_n)\in \Omega^n$ and degrees $d = (d_1,\dotsc,d_n)\in\mathbb{Z}^n$, the function $W_{d}(a)$ is defined by 
\begin{equation}\label{ren_def}
W_{d}(a) =   - \pi \sum_{i\neq j} d_i d_j \log \abs{a_i - a_j} + \pi \sum_{i,j} d_i d_j S_\Omega(a_i,a_j).
\end{equation}
Here $S_\Omega \in C^1(\Omega \times \Omega)$ is the regularization of the Green's function for the London equation on $\Omega$, i.e. 
\begin{equation}\label{S_def}
 S_\Omega(x,y) = G(x,y) + \log\abs{x-y},
\end{equation}
where
\begin{equation}
 \begin{cases}
  -\Delta_x G(x,y) + G(x,y) = 2\pi \delta_{y} & \text{in } \Omega \\
  G(x,y) =0 & \text{for } x\in \partial \Omega, y\in\Omega.
 \end{cases}
\end{equation}
See \cite{ss_book} for further discussion of the magnetic renormalized energy and for proof of \eqref{en_split_intro}.  

Throughout the paper we will use the notation $o(1)$ to refer to any quantity that vanishes as $\ep \rightarrow 0$.  Similarly, $O(1)$ refers to a quantity that stays bounded.  For two quantities $a_\ep$, $b_\ep$, we will employ the notation $a_\ep \ll b_\ep$ to mean that $a_\ep/b_\ep = o(1)$, and we write $a_\ep \asymp b_\ep$  if $a_\ep / b_\ep = O(1)$ and $b_\ep/a_\ep=O(1)$. We will also employ the standard convention of using the letter $C$ to denote a generic positive constant that may change from from line to line.

%%%%%%%%%%%%%%%%%%%%%%%%%%%%%%%%%%%%%%%%%%%%%%%%%%%%%%%%%%%%%%%%%%%%%%%%%%%%%
\subsection{Known results and expectations}
%%%%%%%%%%%%%%%%%%%%%%%%%%%%%%%%%%%%%%%%%%%%%%%%%%%%%%%%%%%%%%%%%%%%%%%%%%%%%

The main interest in studying \eqref{tdgl} is to derive the dynamics of the vortices associated to sequences of solutions $(u_\ep,A_\ep,\Phi_\ep)$ in the limit as $\ep \rightarrow 0$.  In the non-magnetic case, this was accomplished   under various assumptions and for varying lengths of time \cite{bcps,rub_stern,lin_1,js_2,bos_col,bos_conv,bos,serf_1,serf_2}.  The magnetic case with $I_{ex}=J_{ex}=0$ and $H_{ex} = h_{ex}(\ep)$ a constant depending on $\ep$ was studied in \cite{spirn} for $h_{ex}(\ep)$ fixed and in \cite{ss_gamma} for $h_{ex}(\ep) = \beta \ale$ with $\beta >0$.  In this setting the equations \eqref{tdgl} constitute the $L^2$ gradient flow of the full Ginzburg-Landau energy, 
\begin{equation}\label{full_energy}
 G_\ep(u,A) = \hal \int_\Omega \abs{\nab_A u}^2 + \abs{h - h_{ex}(\ep)}^2 +\frac{(1-\abs{u}^2)^2}{2\ep^2}.
\end{equation}
This leads to energy dissipation, $G_\ep(u_\ep(t),A_\ep(t)) \le G_\ep(u_\ep(0),A_\ep(0))$ for all $t\ge 0$, which provides the technical advantage of a priori control of the energy and the number of vortices.  

The general scheme found in these papers is that the vortices do not move until after an amount of time of order $\lep = \ale/h_{ex}(\ep)$.  When $h_{ex}(\ep)$ is fixed, this means that in the limit $\ep \rightarrow 0$ the vortices cannot move at all since $\lep \rightarrow \infty$.  In this case $(u_\ep,A_\ep,\Phi_\ep)\rightarrow (u_*,A_*,\Phi_*)$ in some sense, and it is possible to pass to the limit in \eqref{tdgl} to derive the dynamics for $(u_*,A_*,\Phi_*)$.  To see vortex motion, the solutions are accelerated in time at scale $\lep$ according to 
\begin{equation}
 u_\ep(x,t) \mapsto u_\ep(x,\lep t), A_\ep(x,t) \mapsto A_\ep(x,\lep t), \Phi_\ep(x,t) \mapsto \Phi_\ep(x,\lep t).
\end{equation}
In this scaling, the vortices do move in the limit, and their dynamics are governed by  
\begin{equation}
 \dot{a}_i(t) = -\frac{1}{\pi} \nab_{a_i} W_d(a(t)) - 2 d_i h_{ex} \nab H_0(a_i(t)),
\end{equation}
where $a_i(t)\in \Omega$ is the location of the $i^{th}$ vortex, $d_i\in \mathbb{Z}$ is its degree, $\nab_{a_i} W_d$ is the derivative of magnetic renormalized energy \eqref{ren_def} with respect to $a_i \in \Omega$, and $\nab H_0$ is a purely magnetic forcing term with $H_0$ the solution to the London equation
\begin{equation}
 \begin{cases}
  -\Delta H_0 + H_0 = 0 &\text{in }\Omega \\
   H_0 = 1 & \text{on }\partial \Omega.
 \end{cases}
\end{equation}
As such, the limiting dynamics are a gradient flow of the energy
\begin{equation}
 W_{d,h_{ex}}(a) := W_d(a) + 2\pi h_{ex} \sum_{i=1}^n d_i H_0(a_i),
\end{equation}
the latter term of which is the interaction energy between the vortices and the applied magnetic field.  

In \cite{ss_gamma}, the choice $h_{ex} = \beta \ale$ implies $\lep = 1/\beta$, which allows the vortices to move in the original time scale.  The interpretation of this is that the interaction energy between the vortices and the magnetic field  is sufficiently strong to induce vortex motion in the original time scale.  The resulting motion corresponds to the gradient flow of this interaction energy, and is written $\dot{a}_i = -2\beta d_i \nab H_0(a_i)$.  Note that in this case the vortices do not interact with each other in the sense that the motion of the point $a_i$ does not depend on the points $a_j$ for $j\neq i$.

In the case $J_{ex} \neq 0$ or $I_{ex}\neq 0$ the gradient flow structure of the equations \eqref{tdgl} breaks down.  Energy does not dissipate,  and we can no longer expect the limiting dynamics of the vortices to be a gradient flow.  This creates serious difficulties  in applying the standard Ginzburg-Landau toolboxes, which rely crucially on precise knowledge of the energy.  In fact, the applied boundary current is expected to introduce two novel features to the dynamics, both of which have been observed in the numerical simulations of \cite{du_2,du_3,du_gray}. First, the applied current generates an electric field in $\Omega$, and the vortices feel a Lorentz force perpendicular to this field.  Second, and more drastic, a sufficiently strong applied current is expected to create and destroy vortices near the boundary, thereby injecting or removing large amounts of energy from the system.

%%%%%%%%%%%%%%%%%%%%%%%%%%%%%%%%%%%%%%%%%%%%%%%%%%%%%%%%%%%%%%%%%%%%%%%%%%%%%
\subsection{Summary of main results }
%%%%%%%%%%%%%%%%%%%%%%%%%%%%%%%%%%%%%%%%%%%%%%%%%%%%%%%%%%%%%%%%%%%%%%%%%%%%%

In this paper we analyze sequences of solutions $(u_\ep,A_\ep,\Phi_\ep)$ to \eqref{tdgl} as $\ep \rightarrow 0$ in both the original and accelerated time scales.  Our aim is to show that the applied boundary currents $I_{ex}$ and $J_{ex}$ induce  Lorentz forcing terms in the limiting vortex dynamics for the accelerated time scale, and to identify the structure of the Lorentz forces.   We make the structural assumption that
\begin{equation}
 J_{ex} = j_{ex} J,  H_{ex} = h_{ex} H, \text{ and } I_{ex} = h_{ex} I
\end{equation}
for field strengths $j_{ex} = j_{ex}(\ep)\ge 0$ and $h_{ex}=h_{ex}(\ep) \ge 0$, $J:\partial \Omega \rightarrow \Rn{2}$ a smooth, fixed vector field, and $H:\partial \Omega \rightarrow \Rn{}$ the smooth trace onto $\partial \Omega$ of the solution to the static exterior Maxwell equations $\nb H = -I$.  When $I =0$ we will assume that $H =1$, corresponding to a uniform applied magnetic field.  For notational convenience we will not write the $\ep$ dependence for the parameters $j_{ex}$ or $h_{ex}$.  We shall consider four distinct regimes for the parameters:
\begin{equation}\label{regimes}
\begin{tabular}{ll} 
Regime 1: &  $h_{ex} = j_{ex} =1 $    \\ 
Regime 2: &   $0 \le h_{ex} \ll j_{ex} \ll \ale^{1/9}$                  \\ 
Regime 3: & $0 \le j_{ex} \ll h_{ex} \ll \ale^{1/9}$       \\ 
Regime 4: & $1 \ll h_{ex} \asymp j_{ex} \ll\ale^{1/9}$.
\end{tabular} 
\end{equation}
The first regime handles the choice of any $j_{ex}$ and $h_{ex}$ fixed with respect to $\ep$  since  we may simply rescale $J$ and $I$ to set $j_{ex}=h_{ex}=1$.  In the second two cases, at least one parameter blows up, but one dominates the other.  In the fourth case both blow up but are of the same order.  The upper bound by $\ale^{1/9}$ is purely technical, being required in the proof of the dynamical law in the accelerated time scale.  We define the dominant field strength via
\begin{equation}\label{k_def}
 k_{ex} := \max\{j_{ex},h_{ex}\}.
\end{equation}

The main thrust of the paper is to deal with the complications caused by the dynamics no longer being dissipative or even conservative (as with the Schr\"{o}dinger flows associated to $F_\ep$).  We show that, while the energy does not necessarily decrease, it cannot increase too quickly.  This allows us to identify a time scale depending on $k_{ex}$ in which the number of vortices remains constant, but the vortices move in $\Omega$, exhibiting the additional Lorentz force drifts due to the applied currents.  The control of the energy growth is far from trivial; indeed, we see that standard tricks (e.g. Gronwall) are insufficient for getting precise estimates of the energy.  It is only via the identification of some very special structure in the current forcing terms that we are able to get the delicate estimates required. 

To understand how the free energy of solutions evolves in time, we introduce a splitting of the solutions into a topological (i.e. generated by vortices) component and an applied current component.  This is similar to a technique employed in \cite{bos_col,bos_conv,bos} for the non-gauged case.  To motivate the splitting, we take the curl of the second equation in \eqref{tdgl} to see that 
\begin{equation}
\begin{cases}
 \dt h_\ep - \Delta h_\ep + h_\ep = \mu(u_\ep,A_\ep)  & \text{in }\Omega \\
 h_\ep= h_{ex} H & \text{on } \partial\Omega.
\end{cases}
\end{equation}
We can then split $h_\ep$ according to $h_\ep = h_{ex}h_0 + h_\ep'$ where 
\begin{equation}\label{h0_def}
 \begin{cases}
   -\Delta h_0 +h_0= 0 & \text{in } \Omega \\
   h_0 =H & \text{on } \partial \Omega.
 \end{cases}
\end{equation}
Then $h_0$ is the static contribution of the applied magnetic field, and $h_\ep'$, which satisfies  $h_\ep' =0$ on $\partial \Omega$, is the dynamic part of the induced magnetic field generated by the vortices.   Notice that 
\begin{equation}
 h_\ep' = h_\ep - h_{ex} h_0 = \curl{A_\ep} - h_{ex} \Delta h_0 = \curl{A_\ep} - h_{ex}\curl{\nb h_0} = \curl(A_\ep- h_{ex} \nb h_0).
\end{equation}
This suggests defining the modified vector potential $B_\ep:= A_\ep - h_{ex} \nb h_0$, which we expect to be the part of the vector potential generated by vortices.  Then on $\partial \Omega$
\begin{equation}
 \nab_{B_\ep} u_\ep \cdot \nu = \nab_{A_\ep} u_\ep \cdot \nu + i u_\ep  h_{ex}\nb h_0 \cdot \nu = iu_\ep J_{ex}  \cdot \nu + i u_\ep  \nb H_{ex} \cdot \nu = iu_\ep(J_{ex} - I_{ex})\cdot \nu,
\end{equation}
which shows that $J_{ex}$ acts as a sort of current.  It is then useful to modify $u_\ep$ in such a way to turn the inhomogeneous Neumann boundary condition into a homogeneous one.  To do so we define $f_1$ and $f_0$ as the solutions to 
\begin{equation}\label{fi_def}
 \begin{cases}
   -\Delta f_1 +f_1= -\Delta f_0 +f_0=0 & \text{in } \Omega \\
   \nab f_1 \cdot \nu = J \cdot \nu & \text{on } \partial \Omega \\
   \nab f_0 \cdot \nu = I \cdot \nu & \text{on } \partial \Omega
 \end{cases}
\end{equation}
and write 
\begin{equation}\label{f_def}
f := j_{ex} f_1 - h_{ex} f_0.
\end{equation}
Then the modified order parameter $v_\ep := u_\ep e^{-if}$ satisfies the homogeneous boundary condition
\begin{equation}\label{hom_swap}
 \nab_{B_\ep} v_\ep \cdot \nu = e^{-if} \nab_{B_\ep} u_\ep \cdot \nu - i v_\ep \nab f \cdot \nu = iv_\ep (J_{ex} - I_{ex} - \nab f)\cdot \nu =0.
\end{equation}
It will be convenient to introduce the forcing vector field
\begin{equation}\label{Z_def}
 Z_\ep := j_{ex} \nab f_1 - h_{ex} \nab f_0 - h_{ex} \nb h_0.
\end{equation}

In studying $v_\ep,B_\ep$, we are led by the equations \eqref{tdgl} to consider the evolution of a modification of the standard Ginzburg-Landau free energy $F_\ep$:
\begin{equation}
\tilde{F}_\ep(v_\ep,B_\ep) := F_\ep(v_\ep,B_\ep) +  \frac{1}{2} \int_\Omega \abs{v_\ep}^2  \abs{Z_\ep}^2.
\end{equation}
The study of $\tilde{F}_\ep(v_\ep,B_\ep)$ points to a natural choice of gauge: one in which $\Phi_\ep = f,$ which we refer to as the $\Phi=f$ gauge.  

Many of the arguments in this paper rely crucially on the initial data satisfying a well-preparedness condition on the energy.  In particular, for a sequence $(v_\ep,B_\ep)$ and a constant $C_0>0$, we say that initial data $(v_\ep(0),B_\ep(0))$ are well-prepared at order $C_0$ if $\mu(v_\ep(0),B_\ep(0)) \rightarrow 2\pi \sum_{i=1}^n d_i \delta_{a_i}$ with $d_i =\pm 1$ and 
\begin{equation}\label{well_prepared_def}
 \tilde{F}_\ep(v_\ep(0),B_\ep(0)) \le \pi n \ale +  W_{d}(a) + n\gamma + \frac{1}{2} \int_\Omega \abs{Z_\ep}^2 +  C_0 k_{ex},
\end{equation}
where  $W_{d}(a)$ is the renormalized energy defined by \eqref{ren_def}, $Z_\ep$ is given by \eqref{Z_def}, $k_{ex}$ is defined by \eqref{k_def}, and $\gamma$ is a fixed constant (see Lemma \ref{lower_bound}).  We note that by adapting results in the literature (\cite{ss_book} for example) we may construct initial data satisfying these hypotheses.

The utility of studying $v_\ep, B_\ep$ rather than $u_\ep, A_\ep$ lies in a novel observation on the structure of a term arising in the equation for the evolution of the modified energy.  We find that if $(u_\ep,A_\ep,\Phi_\ep)$ solve \eqref{tdgl} in the $\Phi=f$ gauge, then for $v_\ep=u_\ep e^{-i f}$, $B_\ep=A_\ep - h_{ex}\nb h_0$, 
\begin{equation}
 \dt \tilde{F}_\ep(v_\ep,B_\ep) + \int_\Omega \abs{\dt v_\ep}^2 + \abs{\dt B_\ep}^2 = \int_\Omega V_\ep \cdot Z_\ep, 
\end{equation}
where $V_\ep$ is the ``velocity component'' of the full space-time Jacobian associated to $(v_\ep,B_\ep)$ (see Section \ref{energy_evolution_section}).  The extra structure of the interaction term, $V_\ep \cdot Z_\ep$, is the key to controlling the growth of the modified energy because of  estimates for $V_\ep$ proved in \cite{ss_prod} (recorded here in Proposition \ref{prod_est}).  Using these estimates, we can prove that if the initial data is well-prepared at order $C_0$, then in an amount of time of order   $\lep := \frac{\ale}{k_{ex}},$ the modified energy $\tilde{F}_\ep(v_\ep,B_\ep)$ can increase at most by an amount $2C_0 k_{ex} $.  

\begin{thm_intro}[proved later as Theorem \ref{time_bound}]
Suppose that the initial data $(v_\ep(0),B_\ep(0))$ are well-prepared at order $C_0$, as defined by \eqref{well_prepared_def}.  Let $\lep = \ale/k_{ex}$.  Then there exists a constant $T_0>0$ so that, as $\ep \rightarrow 0$, 
\begin{equation}
\tilde{F}_\ep(v_\ep,B_\ep)(t) \le \tilde{F}_\ep(v_\ep,B_\ep)(0) + 2C_0 k_{ex} \text{  for all } t\in [0,T_0 \lep ], \text{ and}
\end{equation}
\begin{equation}
 \int_0^{T_0 \lep}  \int_\Omega \abs{\dt v_\ep}^2 + \abs{\dt B_\ep}^2   \le 2C_0 k_{ex} .
\end{equation}
\end{thm_intro}

When $j_{ex}=h_{ex}=1$ we can derive the limiting dynamics in the original time scale for the pair $(v_\ep,B_\ep)$ in essentially the same manner as in \cite{spirn}.  %The bounds of the previous theorem can be used to show directly that the vortices do not move.  Then we apply energy localization results, which show that the energy in balls centered at the vortex locations is of order $\pi n \ale$.  This yields bounds (independent of $\ep$) for $(v_\ep,B_\ep)$ in certain Sobolev spaces, which are then used to pass to the limit.  
Since $u_\ep = v_\ep e^{if}$ and $A_\ep = B_\ep + \nb h_0$, we then immediately get the limiting dynamics for $(u_\ep,A_\ep)$ as well.  

\begin{thm_intro}[proved later as Theorem \ref{original_dynamics}]\label{o_d_intro}
Let $j_{ex}=h_{ex} =1$ (parameter regime $1$).  Suppose that the initial data are well-prepared at order $C_0$ as defined by \eqref{well_prepared_def}.  Then on any fixed time interval $[0,T]$ the following hold.
 \begin{enumerate}
 \item The vortex locations  do not move in time, i.e. $a_i(t) = a_i(0)=a_i$ for $t\in[0,T]$.
 
\item $u_\ep \rightharpoonup u_*$ weakly in $H^1_{loc}(\Omega \backslash \{a_i\} \times [0,T])$, where 
\begin{equation}
 u_* = \prod_{i=1}^n \left(\frac{x-a_i}{\abs{x-a_i}}\right)^{d_i} e^{i\psi_* + i f} := e^{i\Theta_a + i \psi_*+ if},
\end{equation}
$f$ is defined by \eqref{f_def}, and $\psi_*$ is a single-valued function on $\Omega\times [0,T]$ satisfying 
\begin{equation}
\begin{cases}
 \dt \psi_* - \Delta \psi_* + \psi_* = \psi_*(0)& \text{in } \Omega \\
 \nab \psi_* \cdot \nu = -\nab \Theta_a \cdot \nu & \text{on } \partial \Omega.
\end{cases}
\end{equation}

\item $A_\ep \rightarrow A_*$  in $L^2(\Omega \times [0,T])$ and $h_\ep \rightharpoonup h_* = \curl{A_*}$ weakly in $L^2(\Omega \times [0,T]])$.  The function $h_*$ satisfies
\begin{equation}
\begin{cases}
 \dt h_* - \Delta h_* + h_* = 2\pi \sum_{i=1}^n d_i \delta_{a_i} &  \text{in } \Omega \\
 h_* = H_{ex} & \text{on } \partial \Omega.
\end{cases}
\end{equation}
\end{enumerate}
\end{thm_intro}

\begin{remark}
If $J  =I =0$, then $f=0$, and we recover a result from \cite{spirn}.  However, in the second item above, the equation satisfied by  $\psi_*$ is different from the equation satisfied by $\psi_*$ in \cite{spirn}.  The source of this disparity is the difference in choice of gauge.  We work in a gauge where $\Phi=f$, whereas \cite{spirn} utilizes the Lorentz gauge, $\Phi+\diverge{A}=0$.  Formally changing $(u_*,A_*,f)$ to the Lorentz gauge shows that our result is consistent with that of \cite{spirn}.
\end{remark}

When $1\ll k_{ex}\ll \ale^{1/9}$ (regimes $2, 3,$ and $4$), the time scale $\lep$ is much smaller than $\ale$, and the vortices begin to move sooner.  However, the limitation on the size of $k_{ex}$ still means that $\lep\rightarrow \infty$ so that it takes infinitely long for the vortices to begin moving in the limit $\ep \rightarrow 0$.  In this case it is possible to extend the previous theorem to show that the vortices do not move in the limit, but unfortunately, the proofs of the second and third items of the theorem break down when $1\ll k_{ex}$, so we can derive no information on the existence or structure of limits of $v_\ep$, $B_\ep$, or $\curl{B_\ep}$ in the original time scale. 

The main result of the paper considers the solutions accelerated in time at scale $\lep$, i.e. we make the substitutions
\begin{equation}\label{acc_def}
v_\ep(x,t) \mapsto v_\ep(x, \lep t) \text{ and } B_\ep(x,t) \mapsto B_\ep(x, \lep t). 
\end{equation}
In this scaling the vortices move along well-defined, continuous trajectories.  To show that the vortex trajectories are differentiable and to derive the limiting law governing their dynamics, we pass to the limit in a localized version of the evolution equation for the modified energy.  %To do so we first prove various convergence results for the (properly normalized) modified energy density, the stress-energy tensor associated to $(v,B)$, and the full space-time Jacobian measure.  
We show that, as predicted, the applied boundary currents $J_{ex}$ and $I_{ex}$ exert  Lorentz forces on the vortices in addition to the forcing term from the magnetic renormalized energy that was identified in \cite{spirn,ss_gamma}.  The exact form of the limiting law depends on the parameter regime.

\begin{thm_intro}[proved later in Lemma \ref{vortex_path} and Theorem \ref{dynamical_law}]\label{dyn_intro}
$\text{}$\\
Suppose that the initial data are well-prepared at order $C_0$ and that the solutions have been accelerated in time at scale $\lep$ according to \eqref{acc_def}.  Suppose further than the initial vortex locations are separated from each other and the boundary by a distance at least $\sigma_0>0$.  Then for $0<\sigma_* < \sigma_0$ there exists a time $T_* = T_*(\sigma_*)\in(0,T_0]$ 
and $n$ differentiable functions $a_i:[0,T_*]\rightarrow \Omega$ satisfying the following.
\begin{enumerate}
 \item For each time $t\in[0,T_*]$ there is a degree $d_i$ vortex located at $a_i(t)$, i.e. the $n$ initial vortices move along the trajectories $a_i$.
 \item The vortices are separated from each other and the boundary  by a distance at least $\sigma_*$ for all time $t\in[0,T_*]$.  In other words, the time $T_*$ is chosen to precede the first time at which a collision occurs or a vortex meets the boundary. 
 \item If $h_{ex} = j_{ex} =1$, then $\lep = \ale$ and the trajectories satisfy the dynamical law
\begin{equation}
 \dot{a}_i(t) = -\frac{1}{\pi} \nab_{a_i}W_{d}(a(t)) -2d_i (\nab h_0(a_i(t)) - \nab^\bot f_0(a_i(t)) +  \nab^\bot f_1 (a_i(t)))
\end{equation}

 \item If $0 \le h_{ex} \ll j_{ex} \ll \ale^{1/9}$, then $\lep = \ale/j_{ex}$ and the trajectories satisfy the dynamical law
\begin{equation}
 \dot{a}_i(t) = -2 d_i \nab^\bot f_1(a_i(t)).
\end{equation}

 \item If $0 \le j_{ex} \ll h_{ex} \ll \ale^{1/9}$, then $\lep = \ale/h_{ex}$ and the trajectories satisfy the dynamical law 
\begin{equation}
 \dot{a}_i(t) = -2 d_i \nab h_0(a_i(t)) + 2d_i \nab^\bot f_0(a_i(t)).
\end{equation}

 \item If $1 \ll h_{ex} \asymp j_{ex} \ll\ale^{1/9}$, $ j_{ex}/k_{ex} \rightarrow \alpha$, and $h_{ex}/k_{ex} \rightarrow \beta$, then $\lep = \ale/k_{ex}$ and the trajectories satisfy the dynamical law
\begin{equation} 
 \dot{a}_i(t) = -2 d_i  \alpha \nab^\bot f_1(a_i(t))  -2d_i \beta  \left(\nab h_0(a_i(t)) -\nb f_0(a_i(t))   \right).
\end{equation}
\end{enumerate}

\end{thm_intro}

\begin{remark}
$\;$
\begin{enumerate}

\item Since $f$ is completely determined by the applied boundary currents $J_{ex}$ and $I_{ex}$, we see that the novel forcing terms in the dynamics, $2 d_i  \nab^\bot f_k (a_i(t))$, $k=0,1$ are really forces induced by the applied currents.  In this way we identify the new terms with the predicted Lorentz forces.  The Lorentz force induced by the applied normal current, $I_{ex}$, is always accompanied by a corresponding magnetic force $2 d_i \nab h_0$.  On the other hand, the Lorentz force from $J_{ex}$ comes unaccompanied by a magnetic force, which indicates a fundamental difference between the currents $J_{ex}$ and $I_{ex}$. 

\item In the case $1\ll k_{ex}$ the magnetic renormalized energy disappears in the limiting dynamical law, and the motion of any given vortex is independent of the motion of the others.  This is  consistent with what was mentioned above about the form of the limiting dynamics for the gradient flow when $h_{ex} = O(1)$ as compared to  when $h_{ex} = \beta \ale$. 

\item In the third parameter regime, the current $J_{ex}$ dominates the dynamics and the limiting law is a perpendicular gradient flow.  This implies that the dynamics of the limiting vortices  possess $n$ conserved quantities: $ \mathcal{H}_i := f_1(a_i) \text{ for }i=1,\dotsc,n.$ This flow keeps the vortices confined to the level sets of the function $f_1$.  Since we do not require $\int_{\partial \Omega} J \cdot \nu=0$ (if $J\cdot \nu$ models a surface charge this should not hold), it is possible to choose $J$ so that the level sets of $f_1$ form closed curves in $\Omega$, which  gives rise to periodic motion of the vortices if no collisions occur. 

\item Working on the time interval $[0,T_*]$ prevents the vortices from colliding with each other or the boundary, and by letting $\sigma_*$ tend to $0$, we can derive the dynamics up to the first collision time.   Unfortunately, our techniques break down at a collision, and we can say nothing about what happens after.  If a collision between two vortices occurs, it may be possible to use a deeper, more refined analysis as in \cite{bos_col,bos_conv,bos,serf_2} to understand the dynamics afterward.  On the other hand, if a vortex collides with the boundary, nothing in our analysis excludes the possibility of that vortex disappearing and another one nucleating somewhere else on $\partial \Omega$. 

\item Setting $j_{ex}=1$, $J=I=0$, and $H = 1$, the second item of this theorem partially bridges the gap between what is known about the limiting dynamical law for the gradient flow dynamics when $h_{ex}=O(1)$ \cite{spirn} and when $h_{ex} =\beta \ale$ \cite{ss_gamma}.  

\end{enumerate}
\end{remark}

\subsection{Plan of the paper}
%%%%%%%%%%%%%%%%%%%%%%%%%%%%%%%%%%%%%%%%%%%%%%%%%%%%%%%%%%%%%%%%%%%%%%%%%%%%%

The paper is organized as follows. In Section \ref{energy_evolution_section} we analyze the dynamics of the modified energy, proving it does not increase too quickly.  In Section \ref{original_time_section} we derive the limiting dynamics in the original time scale. In Section \ref{accelerated_time_section} we study the accelerated solutions and derive the limiting dynamics for the vortices, identifying the force induced by the applied boundary currents $J_{ex}$ and $I_{ex}$.  

At the end of the paper we present two appendices.  Appendix \ref{well_posed_section} contains a few remarks on the well-posedness and regularity of the equations \eqref{tdgl} with the new boundary condition $\nab_A u \cdot \nu = iu J_{ex}\cdot \nu$ as well as some a priori estimates needed in a couple places in the paper.  We will also need several results from the static analysis of the Ginzburg-Landau energy.  These are collected in Appendix \ref{static_section}.    

\vskip .5cm

{\it Acknowledgments}: Part of this work was completed as a component of my Ph.D. thesis at the Courant Institute.  I would like to thank my advisor, Sylvia Serfaty.  Thanks also to Fang-Hua Lin and Qiang Du for their helpful suggestions.   Finally, I would like to thank the NSF for the support of the GRF during my Ph.D. studies, and also the Laboratoire Jacques-Louis Lions at Paris 6 for their hospitality during my visit.

%%%%%%%%%%%%%%%%%%%%%%%%%%%%%%%%%%%%%%%%%%%%%%%%%5555
\section{Energy evolution}\label{energy_evolution_section}
%%%%%%%%%%%%%%%%%%%%%%%%%%%%%%%%%%%%%%%%%%%%%%%%%%%%%

%%%%%%%%%%%%%%%%%%%%%%%%%%%%%%%%%%%%%%%%%%%%%%%%%%%%%
\subsection{Evolution of the modified energy}
%%%%%%%%%%%%%%%%%%%%%%%%%%%%%%%%%%%%%%%%%%%%%%%%%%%%%

The obvious starting point for an analysis of the behavior of a sequence of solutions $(u_\ep,A_\ep,\Phi_\ep)$ to \eqref{tdgl} as $\ep \rightarrow 0$ is an examination of how the energy $F_\ep(u,A)$ evolves in time.  To understand this, we consider the Ginzburg-Landau free energy density
\begin{equation}
 g^0_\ep(u,A): = \hal \left(\abs{\nab_A u}^2 + \frac{1}{2\ep^2}(1-\abs{u}^2)^2 + \abs{\curl{A} }^2  \right).
\end{equation}
Then $\dt F_\ep(u,A) = \int_\Omega \dt g^0_\ep(u,A)$, and we are led to consider the time derivative of the energy density.

\begin{lem}\label{en_evolve}
 For any triple $(w,B,\xi)$ (not necessarily solutions), we have that 
\begin{equation}
\begin{split}
 \dt g^0_\ep(w,B) & = \diverge(\partial_\xi w,\nab_B w) + \curl((\curl{B})(\dt B + \nab \xi) )  \\
 & - \left(\partial_\xi w , \Delta_B w + \frac{w}{\ep^2}(1-\abs{w}^2)\right) - (\dt B + \nab \xi) \cdot (\nab^\bot \curl{B} + (iw,\nab_B w)).
\end{split}
\end{equation}
\end{lem}
\begin{proof}
The result follows from a simple calculation and an application of the following lemma, which records the commutator relations for the covariant derivatives.
\end{proof}

\begin{lem}\label{commutator}
For any triple $(w,B,\xi)$ (not necessarily solutions), the following commutation relations hold for the covariant derivatives, $\partial_j^B := \partial_j - iB_j$, $\partial_\xi := \dt  + i \xi$:
\begin{equation}
\begin{split}
& \partial_2^B \partial_1^B w - \partial_1^B \partial_2^B w = i w \curl{B} \\
& \nab_B  \partial_{\xi} w - \partial_\xi \nab_B w = i w (\dt B + \nab \xi) .
\end{split}
\end{equation}
\end{lem}

So, to understand $\dt F_\ep(u,A)$, we apply Lemma \ref{en_evolve} to $(u,A,\Phi)$ and integrate over $\Omega$.  The integrals of the $\diverge$ and $\curl$ terms yield a non-vanishing boundary integral because of the conditions $\nab_{A} u \cdot \nu = i u J_{ex} \cdot \nu$  and $h= H_{ex}$.  In particular, we find that for solutions to \eqref{tdgl}
\begin{equation}
 \dt F_\ep(u,A) + \int_\Omega \abs{\partial_{\Phi} u}^2 + \abs{E}^2 = \int_{\partial \Omega} (iu,\partial_{\Phi} u) J_{ex} \cdot \nu - H_{ex} E\cdot \tau.
\end{equation}
This boundary integral is inconvenient to work with, so we pursue an alternate strategy for studying the energy evolution.  The idea is to modify  $u$ and $A$ in a manner that turns the inhomogeneous boundary condition $\nab_{A} u \cdot \nu = i u J_{ex} \cdot \nu$ and $h = H_{ex}$ into  homogeneous ones.  This is accomplished by introducing the modified order parameter $v = u e^{-i f}$ and the modified vector potential $B = A - h_{ex} \nab^\bot h_0$, where $f$ and $h_0$ are the functions defined by \eqref{f_def} and \eqref{h0_def} respectively.  We emphasize that we are \emph{not} making a gauge change since we \emph{do not} make either of the changes $A \mapsto A  +\nab f$ or $\Phi \mapsto \Phi - \dt f$.  The choice of $f$ and $h_0$ then imply that 
 $\nab_B v \cdot \nu = 0$ and $\curl{B} =0$.  Throughout the rest of the paper we will write
\begin{equation}\label{h_prime_def}
 h' := \curl{B} = \curl(A - h_{ex} \nb h_0) = h - h_{ex} h_0.
\end{equation}
The trade-off for switching to homogeneous boundary conditions is that the equations pick up forcing terms.  We record the equation satisfied by $(v,B,\Phi)$ now.

\begin{lem}\label{eqn_swap}
 Suppose the triple $(u,A,\Phi)$ solve \eqref{tdgl}.  Let $v = u e^{-if}$, $B= A - h_{ex}\nab^\bot h_0$.  Then the triple $(v,B,\Phi)$ solve 
\begin{equation}\label{mod_eqn_1}
 \dt v + i \Phi v = \Delta_B v + \frac{v}{\varepsilon^2}(1-\abs{v}^2) 
+ 2i\nab_B v \cdot Z_\ep + iv f - v  \abs{Z_\ep  }^2
\end{equation}
and
\begin{equation}
 \dt B + \nab \Phi = \nab^\bot h' + (iv,\nab_B v) +\abs{v}^2 \nab f   - (\abs{v}^2 -1) h_{ex} \nab^\bot h_0
\end{equation}
along with the homogeneous boundary conditions
\begin{equation}\label{mod_eqn_2}
 \begin{cases}
  \nabla_B v \cdot \nu = 0 &\text{on } \partial \Omega \times \Rn{+} \\
  h' = 0 & \text{on } \partial \Omega \times \Rn{+}. \\
 \end{cases}
\end{equation}
Here we have written $h' = \curl{B}$ and $Z_\ep$ is defined by \eqref{Z_def}.
\end{lem}
\begin{proof}
 A direct calculation gives the PDEs.  The boundary conditions follow via \eqref{hom_swap}.
\end{proof}

In considering the evolution of the energy (in particular, when working with $v,B$) it is most convenient to work in a gauge for which the electric potential $\Phi$ is set equal to the function $f = j_{ex} f_1 -h_{ex} f_0$. 

\begin{lem}\label{f_gauge}
 We can fix a gauge so that $\Phi =  f$ and $B(0)$ satisfies the Coulomb gauge, i.e $\diverge{B(0)}=0$ in $\Omega$ and  $B(0)\cdot \nu =0$ on $\partial \Omega$.
\end{lem}
\begin{proof}
 We make a gauge change $v \mapsto ve^{i\xi}$, $B \mapsto   B + \nab \xi$, $\Phi \mapsto \Psi: = \Phi - \dt \xi$, where $\xi$ is given by
\begin{equation}
 \xi(x,t) = \int_0^t \Phi(x,s) ds - t f(x) + \eta(x)
\end{equation}
and $\eta$ is the solution to 
\begin{equation}
 \begin{cases}
  -\Delta \eta = \diverge{B(0)} & \text{in } \Omega \\
   \nab \eta \cdot \nu = -B(0) \cdot \nu & \text{on } \partial \Omega.
 \end{cases}
\end{equation}
Then $\Psi =  f$ and $(B+\nab \xi)(0)$ satisfies the Coulomb gauge as desired.
\end{proof}

\begin{remark}
In all of what follows we work exclusively in this gauge, which we call the $\Phi=f$ gauge.  We will also cease to refer to solution triples $(v,B,\Phi)$ and instead refer to just $(v,B)$ since $\Phi = f$.
\end{remark}

Now we record the energy evolution equation for $(v,B)$ in the $\Phi=f$ gauge.

\begin{lem}\label{en_evolve_gaugefixed}
Let $(v,B)$ solve \eqref{mod_eqn_1}--\eqref{mod_eqn_2} in the $\Phi=f$ gauge.  Then
 \begin{multline}
  \dt g^0_\ep(v,B) = \diverge(\dt v, \nab_B v) + \curl( h' \dt A) - \abs{\dt v}^2 - \abs{\dt B}^2 \\
 + (\abs{v}^2-1)  Z_\ep \cdot \dt B + 2(\dt v,i\nab_B v) \cdot Z_\ep 
- (v,\dt v)\abs{Z_\ep}^2 
 \end{multline}
\end{lem}
\begin{proof}
 Plug the triple $(v,B,0)$ into Lemma \ref{en_evolve} and then plug in the equations of \eqref{mod_eqn_1}--\eqref{mod_eqn_2} with $\Phi =  f$.
\end{proof}

Since $f$ and $h_0$ do not depend on time, the last term on the right hand side of the last equation is a time derivative.  This leads us to define the modified energy density
\begin{equation}
 \tilde{g}^0_\ep(v,B) := g^0_\ep(v,B) + \frac{\abs{v}^2}{2}  \abs{Z_\ep}^2
\end{equation}
so that $(v,B)$ satisfy
\begin{equation} 
 \dt \int_\Omega \tilde{g}^0_\ep(v,B) + \int_\Omega \abs{\dt v}^2 + \abs{\dt B}^2 = \int_\Omega \left(2(\dt v,i\nab_B v) + (\abs{v}^2-1) \dt B \right)  \cdot   Z_\ep.
\end{equation}
The utility of this equation is that it allows us to identify the two distinct terms that contribute to the change of the modified energy.  We find a standard dissipative term, 
\begin{equation}
 \int_\Omega \abs{\dt v}^2 + \abs{\dt B}^2, 
\end{equation}
which acts to decrease the modified energy.  We also find the interaction term, 
\begin{equation}\label{interaction_def}
 \int_\Omega \left(2(\dt v,i\nab_{B} v) + (\abs{v}^2-1) \dt B \right)  \cdot  Z_\ep,
\end{equation}
which mediates the interaction between $(v,B)$ and the applied boundary currents $J_{ex}, I_{ex}$ and magnetic field $H_{ex}$.  We will eventually see that there is an exceptionally nice structure to the interaction term, but for now we ignore this structure and present a crude preliminary estimate of the growth of $\int_\Omega \tilde{g}^0_\ep(v,B)$.

\begin{prop}\label{crude_energy}
Let $(v,B)$ solve \eqref{mod_eqn_1}--\eqref{mod_eqn_2} in the $\Phi=f$ gauge.  Then for any $t\ge 0$, 
\begin{equation}\label{c_e_0}
 \int_\Omega \tilde{g}^0_\ep(v,B)(t) + \hal \int_0^t \int_\Omega \abs{\dt v}^2 + \abs{\dt B}^2  
\le \exp(C  t) \int_\Omega \tilde{g}^0_\ep(v,B)(0),
\end{equation}
where $C = 4\pnorm{Z_\ep}{\infty}^2$.
\end{prop}

\begin{proof}
Cauchy's inequality implies that 
\begin{equation}\label{c_e_1}
 \int_\Omega 2 (\dt v,i\nab_B v)\cdot Z_\ep  \le \frac{1}{2} \int_\Omega \abs{\dt v}^2 + 4  \pnorm{Z_\ep}{\infty}^2 \int_\Omega \hal \abs{\nab_B v}^2
\end{equation}
and that 
\begin{equation}\label{c_e_2}
\int_\Omega  (\abs{v}^2-1) \dt B \cdot Z_\ep \\
\le \hal \int_\Omega \abs{\dt B}^2 + 2\ep^2  \pnorm{Z_\ep}{\infty}^2 \int_\Omega \frac{(1-\abs{v}^2)^2}{4\ep^2}.
\end{equation}
Plugging \eqref{c_e_1} and  \eqref{c_e_2} into Lemma \ref{mod_en_evolve} then yields the differential inequality 
\begin{equation}
 \dt \int_\Omega \tilde{g}^0_\ep(v,B)  + \hal \int_\Omega \abs{\dt v}^2  +\abs{\dt B}^2 \le 4  \pnorm{Z_\ep}{\infty}^2  \int_\Omega \tilde{g}^0_\ep(v,B).
\end{equation}
An application of Gronwall then proves \eqref{c_e_0}.

\end{proof}

%%%%%%%%%%%%%%%%%%%%%%%%%%%%%%%%%%%%%%%%%%%%%%%%%%%%%
\subsection{Precursors for estimating $B$ in the $\Phi=f$ gauge}
%%%%%%%%%%%%%%%%%%%%%%%%%%%%%%%%%%%%%%%%%%%%%%%%%%%%%

In order to implement a more refined analysis of the interaction term \eqref{interaction_def}, we will require the ability to control the term $\int_\Omega \abs{B}^2$.  Since $\curl{B}$ is gauge invariant and $\diverge{B}$ can be controlled in some gauge choices, this is most naturally accomplished by using a Poincar\'{e} inequality of the form
\begin{equation}
 \int_\Omega \abs{B}^2 \le C \int_\Omega \abs{\diverge{B}}^2 + \abs{\curl{B}}^2.
\end{equation}
Such an inequality fails in general (e.g. $B = \nab h$ for $h$ harmonic), but is available, for instance, in the space $H^1_n(\Omega;\Rn{2}) =  \{B\in H^1(\Omega;\Rn{2})  \;\vert\; B\cdot \nu =0 \text{ on }\partial\Omega\}$, which the magnetic potential $B$ belongs to in both the  Lorentz ($\Phi+\diverge{A}=0$) and Coulomb ($\diverge{A}=0$) gauges.  Unfortunately, in the $\Phi=f$ gauge, it no longer holds that $B \cdot \nu=0$, so we must resort to a version of Poincar\'{e} that also involves a boundary term.  In this section we will record such a Poincar\'{e} inequality, then prove two lemmas needed to conveniently use this inequality in the $\Phi=f$ gauge.

We begin with the aforementioned version of the Poincar\'{e} inequality, as well as another useful estimate.

\begin{prop}\label{boundary_poincare}
Let
\begin{equation}
\mathcal{X}:= \{B\in L^2(\Omega;\Rn{2}) \;\vert\; \curl{B}, \diverge{B} \in L^2(\Omega), \text{ and } B\cdot \nu \in L^2(\partial \Omega)  \}.
\end{equation}
Then there exists a constant $C>0$ so that for any $B\in \mathcal{X}$ it holds that 
\begin{equation}\label{bnd_p_0}
\int_\Omega \abs{B}^2 \le C \left(\int_\Omega \abs{\curl{B}}^2 + \abs{\diverge{B}}^2 + \int_{\partial\Omega}\abs{B \cdot \nu}^2  \right).
\end{equation}
\end{prop}
\begin{proof}
For any $B \in \mathcal{X}$ we solve the elliptic problem
\begin{equation}
\begin{cases}
 \Delta \phi = \diverge{B} &\text{in }\Omega \\
 \nab \phi \cdot \nu = B \cdot \nu &\text{on }\partial \Omega.
\end{cases}
\end{equation}
Then standard elliptic estimates give
\begin{equation}
 \norm{\phi}_{H^1(\Omega)} \le C\left(\norm{\diverge{B}}_{L^2(\Omega)} + \norm{B\cdot \nu}_{H^{-1/2}(\partial \Omega)}   \right) \le C\left(\norm{\diverge{B}}_{L^2(\Omega)} + \norm{B\cdot \nu}_{L^{2}(\partial \Omega)}   \right).
\end{equation}
We may then define the vector field $\tilde{B} = B - \nab \phi$, which satisfies $\diverge{\tilde{B}}=0$ in $\Omega$ and $\tilde{B}\cdot \nu = 0$ on $\partial \Omega$.   The Poincar\'{e} inequality mentioned at the beginning of the section is applicable to $\tilde{B}$ and yields
\begin{equation}
 \norm{\tilde{B}}_{L^2(\Omega)} \le C \norm{\curl{\tilde{B}}}_{L^2(\Omega)} = C \norm{\curl{B}}_{L^2(\Omega)}.
\end{equation}
Combining these two estimates yields the desired inequality.
\end{proof}

\begin{remark}
The Poincar\'{e} inequality recorded in this proposition is not optimal, but is well-suited for our analysis.  The above proof shows that it could be replaced with
\begin{equation}
 \norm{B}_{L^2(\Omega)} \le C\left(\norm{\diverge{B}}_{L^2(\Omega)}+ \norm{\curl{B}}_{L^2(\Omega)} + \norm{B\cdot \nu}_{H^{-1/2}(\partial \Omega)}   \right).
\end{equation}
 
\end{remark}

In order to apply this version of the Poincar\'{e} inequality, we will need estimates for $\diverge{B}$ and $B \cdot \nu$ in the $\Phi=f$ gauge.  These are presented in the following lemma.

\begin{lem}\label{diverge_control}
Let $(v,B)$ solve \eqref{mod_eqn_1}--\eqref{mod_eqn_2} in the $\Phi=f$ gauge. Then
\begin{equation}\label{div_c_00}
 \pnormspace{B \cdot \nu (t)}{2}{\partial\Omega} \le \int_0^t \pnormspace{(\abs{v}^2-1) j_{ex} \nab f_1 \cdot \nu}{2}{\partial\Omega}
\end{equation}
and
\begin{equation}\label{div_c_0}
 \pnormspace{\diverge{B}(t)}{2}{\Omega} \le \int_0^t \pnormspace{\dt v}{2}{\Omega} + \int_0^t \pnormspace{(\abs{v}^2-1)   f}{2}{\Omega}
\end{equation}
\end{lem}
\begin{proof}
 The second equation in \eqref{mod_eqn_1}--\eqref{mod_eqn_2} reads 
\begin{equation}
\dt B = \nab^\bot h' + (iv,\nab_B v) + (\abs{v}^2-1) Z_\ep. 
\end{equation}
Taking the dot product of this equation with the boundary normal $\nu$ and applying the boundary conditions $h' = 0$ and $\nab_B v \cdot \nu =0$ on $\partial \Omega$ yields
\begin{multline}
 \dt B \cdot \nu = \nab^\bot h' \cdot \nu + (iv,\nab_B v\cdot \nu) + (\abs{v}^2-1) ( \nab f -h_{ex} \nb h_0)\cdot \nu \\
 = (\abs{v}^2-1) (j_{ex} J - h_{ex} I + h_{ex} I) \cdot \nu = (\abs{v}^2-1) j_{ex} \nab f_1 \cdot \nu.
\end{multline}
Then, since $B\cdot \nu(0) =0$, we get 
\begin{equation}
 B \cdot \nu (t) = \int_0^t (\abs{v(s)}^2-1) j_{ex} \nab f_1 \cdot \nu,
\end{equation}
from which \eqref{div_c_00} follows.

Taking the divergence of the second Ginzburg-Landau equation, we find that
\begin{equation}
 \dt \diverge{B} +  \Delta f = \diverge((iv,\nab_B v) + \abs{v}^2 Z_\ep   ).
\end{equation}
On the other hand, since $\Delta f = f$, taking $(iv,\cdot)$ with the first equation in \eqref{mod_eqn_1}--\eqref{mod_eqn_2} yields the equality
\begin{equation}
 (iv,\dt v) + \abs{v}^2  f = \diverge((iv,\nab_B v) + \abs{v}^2 Z_\ep )
\end{equation}
so that
\begin{equation}
 \dt \diverge{B} = (iv,\dt v) + (\abs{v}^2-1)  f.
\end{equation}
Since $\diverge{B}(0)=0$, we get
\begin{equation}
 \diverge{B}(t) = \int_0^t (iv,\dt v) +  \int_0^t (\abs{v}^2-1) f,
\end{equation}
which yields \eqref{div_c_0}.
\end{proof}

The next result provides control of the boundary term $(\abs{v}^2-1)\vert_{\partial \Omega}$, which is necessary for the estimate \eqref{div_c_00} to be useful.

\begin{lem}\label{bndry_conv}
Suppose $v:\Omega\rightarrow \mathbb{C}$ (not necessarily a solution) satisfies $\pnorm{v}{\infty}\le K$ and 
\begin{equation}\label{b_con_0}
 \int_\Omega \hal  \abs{\nab \abs{v}}^2 + \frac{(1- \abs{v}^2)^2}{4\ep^2} \le K \ale.
\end{equation}
Then  
 \begin{equation}\label{b_con_00}
  \int_{\partial \Omega} \abs{1-\abs{v}^2}^2 \le C \sqrt{\ep} \ale
 \end{equation}
for a constant $C>0$, depending on $\Omega$ and $K$.
 \end{lem}
\begin{proof}

Using trace and interpolation theory, we have
\begin{equation}
 \norm{1-\abs{v}^2}^2_{L^2(\partial\Omega)} \le C \norm{1-\abs{v}^2}^2_{H^{3/4}(\Omega)} 
\le C \norm{1-\abs{v}^2}_{H^{1}(\Omega)}^{3/2} \norm{1-\abs{v}^2}_{L^{2}(\Omega)}^{1/2}.
\end{equation}
Since $\abs{v}\le K$, we may bound $\abs{\nab(1-\abs{v}^2)}=\abs{2\abs{v}\nab\abs{v}} \le 2K \abs{\nab \abs{v}}$, so that \eqref{b_con_0} implies 
\begin{equation}
 \norm{1-\abs{v}^2}_{H^{1}(\Omega)} \le \left( \int_\Omega  2K \abs{\nab \abs{v}}^2 + (1-\abs{v}^2)^2 \right)^{1/2} \le \sqrt{4K^2 \ale}.
\end{equation}
On the other hand, \eqref{b_con_0} also implies that 
\begin{equation}
 \norm{1-\abs{v}^2}_{L^{2}(\Omega)} \le \sqrt{4K \ep^2 \ale}
\end{equation}
Then \eqref{b_con_00} follows by combining these three bounds.

\end{proof}

%%%%%%%%%%%%%%%%%%%%%%%%%%%%%%%%%%%%%%%%%%%%%%%%%%%%%
\subsection{Control of the modified energy}
%%%%%%%%%%%%%%%%%%%%%%%%%%%%%%%%%%%%%%%%%%%%%%%%%%%%%

Fortunately, the crude Cauchy-Gronwall combination used in the proof of Proposition \ref{crude_energy} is not the only method available for estimating the change of the modified energy.  The second method relies on a more careful analysis of the interaction term \eqref{interaction_def}.  We begin this line of reasoning by examining the space-time 1-form defined for any triple $(w,B,\xi)$ by
\begin{equation}
 \omega := ((iw,\partial_\xi w) - \xi)dt + ((iw,d_B w) +B),
\end{equation}
where we have identified the vector $B$ with the spatial 1-form $B=B_1 dx_1 + B_2 dx_2$ and written the 1-form $d_B w := d_{sp}w + iwB$ for $d_{sp}$ the spatial exterior derivative in $\Omega$.  This 1-form is a modification of the usual space-time current, $(iw,\partial_\xi w) dt + (iw,d_B w)$.  Taking the exterior derivative of $\omega$, we find the space-time Jacobian measure
\begin{equation}
 \mathcal{J} = d\omega = V_1 dx_1\wedge dt + V_2 dx_2 \wedge dt + \mu dx_1\wedge dx_2,
\end{equation}
where $\mu = \curl((iw,\nab_B w) +B)$ is the usual spatial vorticity measure and the vector 
\begin{equation}\label{V_def}
 V:=(V_1,V_2) =  \nab (iw,\partial_\xi w) - \dt (iw,\nab_B w) - (\dt B + \nab \xi)
\end{equation}
is the ``velocity component'' of the space-time Jacobian measure. We will sometimes write $\mu = \mu(w,B)$ and $V = V(w,B,\xi)$ to emphasize the dependence of $\mu$ and $V$ on $(w,B,\xi)$.  Note, though, that both $\mu$ and $V$ are gauge invariant.

Since $d\circ d=0$, it holds that $d\mathcal{J}=0$, which implies the equation
\begin{equation}
\dt \mu + \curl{V} =0. 
\end{equation}
This relation allows us to identify the vector $V$ with the velocity of the vortices.  To see this, consider the simple example of $\mu(t) = 2\pi \delta_{\gamma(t)}$ for some smooth curve $\gamma:(t_0,t_1)\rightarrow \Omega$.  A straightforward calculation then shows that $V(t) = 2\pi \dot{\gamma}^\bot(t) \delta_{\gamma(t)}$.  Note that even though $V$ actually encodes the perpendicular to $\dot{\gamma}$, we still refer to $V$ as the ``velocity.''

Now we relate the velocity associated to a triple $(w,B,\xi)$ to a quantity that appears in the evolution equation.
\begin{lem}\label{V_rewrite}
 For any triple $(w,B,\xi)$ it holds that
\begin{equation}
 V(w,B,\xi) + \dt B + \nab \xi = 2 (\partial_\xi w,i \nab_B w ) +  \abs{w}^2 (\dt B + \nab \xi).
\end{equation}
\end{lem}
\begin{proof}
 Using Lemma \ref{commutator}, we calculate
\begin{equation}
 \begin{split}
  \nab (iw,\partial_\xi w) & = (\partial_\xi w,i\nab_B w) + (iw, \nab_B \partial_\xi w)   \\
& = (\partial_\xi w,i\nab_B w) + (iw, \partial_\xi  \nab_B w + iw(\dt B + \nab \xi)) \\
& = (\partial_\xi w,i\nab_B w) +  \dt (iw,\nab_B w) - (i\partial_\xi w, \nab_B w) + \abs{w}^2 (\dt B+\nab \xi) \\
& = 2(\partial_\xi w,i\nab_B w) +  \dt (iw,\nab_B w) + \abs{w}^2 (\dt B + \nab \xi).
 \end{split}
\end{equation} 
Hence 
\begin{multline}
 V(w,B,\xi)+(\dt B + \nab \xi) = \nab (iw,\partial_\xi w) - \dt (iw,\nab_B w) \\
= 2 (\partial_\xi w,i\nab_B w) +  \abs{w}^2 (\dt B + \nab \xi). 
\end{multline}

\end{proof}

We now use this lemma  to rewrite the interaction term in the modified energy evolution equation.

\begin{lem}\label{mod_en_evolve}
  Let $(v,B)$ solve \eqref{mod_eqn_1}--\eqref{mod_eqn_2} in the $\Phi=f$ gauge.  Then
\begin{equation}\label{m_e_e_0}
   \dt \tilde{g}^0_\ep(v,B)  = \diverge(\dt v, \nab_B v) + \curl(h' \dt B) 
- \abs{\dt v}^2 - \abs{\dt B}^2 
+ V(v,B,0)\cdot Z_\ep
\end{equation}
and 
\begin{equation}\label{m_e_e_00}
 \dt \int_\Omega \tilde{g}^0_\ep(v,B) + \int_\Omega \abs{\dt v}^2 + \abs{\dt B}^2 = \int_\Omega V(v,B,0) \cdot Z_\ep.
\end{equation}
\end{lem}

\begin{proof}
Use the triple $(v,B,0)$ in Lemma \ref{V_rewrite} to find
\begin{equation}
  (\abs{v}^2-1) \dt B + 2(\dt v,i\nab_B v)  = V(v,B,0). 
\end{equation}
Plugging this into Lemma \ref{en_evolve_gaugefixed} yields \eqref{m_e_e_0}, and \eqref{m_e_e_00} follows by integrating over $\Omega$ and using the boundary conditions.

\end{proof}

The purpose of rewriting the interaction term as $V(v,B,0) \cdot Z_\ep$ is that there are  estimates of the velocity available from \cite{ss_prod}.  These estimates are actually for the velocity vector $V(v,0,0)$, but they will be sufficient for analyzing $V(v,B,0)$.

\begin{prop}[Product Estimate, Theorem 3 of \cite{ss_prod}]\label{prod_est}
 Let $u_\ep:\Omega \times [0,T]\rightarrow \mathbb{C}$ satisfy 
\begin{equation}
\int_\Omega \abs{\nab u_\ep}^2 + \frac{(1-\abs{u_\ep}^2)^2}{2\ep^2}  \le C \ale \text{ for all } t\in[0,T], 
\end{equation}
and
\begin{equation}
 \int_0^T \int_\Omega \abs{\dt u_\ep}^2 \le C \ale.
\end{equation}
Then the following hold.
\begin{enumerate}
 \item There exist $\mu \in L^\infty([0,T];\mathcal{M}(\Omega)) $ and $V=(V_1,V_2)$ with  $V_i \in L^2([0,T];\mathcal{M}(\Omega))$ such that $\dt \mu + \curl{V} = 0$ and $\mu(u_\ep,0) \rightarrow \mu$, $V(u_\ep,0,0) \rightarrow V$ as $\ep \rightarrow 0$ in $(C^{0,\alpha}(\Omega\times[0,T]))^*$ for any $\alpha\in(0,1)$.  Here we have written $\mathcal{M}(\Omega) = (C^0(\Omega))^*$ for the space of bounded Radon measures.

 \item For any $[t_1,t_2] \subseteq [0,T]$ and any $X \in C^0(\Omega \times [t_1,t_2];\Rn{2})$, we have the bound 
\begin{equation}\label{p_e_1}
\hal \abs{\int_{t_1}^{t_2} \int_\Omega V \cdot X} 
\le \liminf_{\ep \rightarrow 0} \frac{1}{\ale}\left( \int_{t_1}^{t_2}  \int_\Omega \abs{\nab u_\ep \cdot X}^2  \int_{t_1}^{t_2}  \int_\Omega \abs{\dt u_\ep}^2  \right)^{1/2}.
\end{equation}

\item The mapping $t \mapsto \langle \mu(t),\xi \rangle$ is in $H^1([0,T])$ for any $\xi \in C_c^1(\Omega)$, and in particular 
\begin{equation}
 \abs{\langle\mu(t_2),\xi\rangle-\langle\mu(t_1),\xi\rangle} 
\le C \sqrt{t_2-t_1} \liminf_{\ep \rightarrow 0} \frac{1}{\sqrt{\ale}} \left( \int_{t_1}^{t_2} \int_\Omega \abs{\dt  u_\ep}^2  \right)^{1/2}
\end{equation}
for any $[t_1,t_2]\subseteq [0,T]$.

\end{enumerate}
\end{prop}

\begin{remark}
In the first item of Proposition \ref{prod_est} the convergence we record is slightly stronger than what is stated in \cite{ss_prod}.  This is valid because the proof in \cite{ss_prod} relies on the method used in \cite{js}, where the stronger convergence result is actually proved.
\end{remark}

With this estimate of the velocity in hand, we can now show more refined control of the growth of the modified energy.  The proof relies critically on the initial data satisfying a well-preparedness assumption given by \eqref{well_prepared_def}.   Recall the dominant field strength, $k_{ex}$, is given by \eqref{k_def}.

\begin{thm}\label{time_bound}
Let $(v_\ep,B_\ep)$ solve \eqref{mod_eqn_1}--\eqref{mod_eqn_2} in the $\Phi=f$ gauge.   Suppose that the initial data $(v_\ep(0),B_\ep(0))$ are well-prepared at order $C_0$, as defined by \eqref{well_prepared_def}.  Define the time scale $\lep := \ale / k_{ex}$.  Then there exists a constant $T_0>0$ such that, as $\ep \rightarrow 0$,
\begin{equation}
 \int_\Omega \tilde{g}^0_\ep(v_\ep,B_\ep)(t) < \int_\Omega \tilde{g}^0_\ep(v_\ep,B_\ep)(0) + 2C_0 k_{ex}  \text{ for all } t\in [0,T_0\lep]
\end{equation}
and 
\begin{equation}
 \int_0^{T_0\lep} \int_\Omega \abs{\dt v_\ep}^2 + \abs{\dt B_\ep}^2  < 2C_0 k_{ex}.
\end{equation}
\end{thm}

\begin{proof}
The proof is inspired by that of Lemma III.1 in \cite{ss_gamma}.  For any $t\ge 0$, consider the two conditions
\begin{equation}\label{cond1}
 \int_\Omega \tilde{g}^0_\ep(v_\ep,B_\ep)(s) < \int_\Omega \tilde{g}^0_\ep(v_\ep,B_\ep)(0) + 2C_0 k_{ex}\text{ for all } s\in [0,t]
\end{equation}
and 
\begin{equation}\label{cond2}
 \int_0^{t} \int_\Omega  \abs{\dt v_\ep}^2 + \abs{\dt B_\ep}^2   < 2C_0 k_{ex}.
\end{equation}
Define 
\begin{equation}
 \alpha_\ep := \sup \{ t\ge 0 \;\vert\; \text{conditions }\eqref{cond1} \text{ and } \eqref{cond2} \text{ hold}   \}.
\end{equation}
Proposition \ref{crude_energy} guarantees the existence of a time $t_\ep>0$ (depending on $\ep$, $j_{ex}$, $h_{ex}$ and $C_0$) such that both conditions hold for $t_\ep$.  Hence $\alpha_\ep > 0$ for each $\ep$.  We will show that actually $\alpha_\ep \ge T_0 \lep$  for some $T_0>0$ as $\ep \rightarrow 0$, thereby proving the theorem.  Suppose now, by way of contradiction, that 
\begin{equation}
\liminf_{\ep\rightarrow 0}  \frac{\alpha_\ep}{\lep} =0.
\end{equation}
We may suppose, up to extraction of a subsequence, that $\alpha_\ep/\lep \rightarrow 0$ as $\ep \rightarrow 0$.

Rescale in time at scale $\alpha_\ep$ by defining $w_\ep(x,t) = v_\ep(x,\alpha_\ep t)$ and $C_\ep(x,t) = B_\ep(x,\alpha_\ep t)$.   By the definition of $\alpha_\ep$, the inequalities
\begin{equation}\label{cond1'}
  \int_\Omega \tilde{g}^0_\ep(w_\ep,C_\ep)(t) \le \int_\Omega \tilde{g}^0_\ep(w_\ep,C_\ep)(0) + 2C_0 k_{ex} \text{ for all } t\in [0,1]
\end{equation}
and 
\begin{equation} \label{cond2'}
 \frac{1}{\alpha_\ep} \int_0^{1} \int_\Omega \abs{\dt w_\ep}^2 + \abs{\dt C_\ep}^2   \le 2C_0 k_{ex}.
\end{equation}
both hold, but at time $t=1$ one of the inequalities must be an equality since $w_\ep$ and $C_\ep$ are smooth.  Our goal is to show that, using the product estimate, neither inequality can fail at time $t=1$, producing the desired contradiction.  In order to be able to apply the product estimate, though, we must first show that its hypotheses are satisfied.  For the rescaled pair $(w_\ep,C_\ep)$, equation \eqref{m_e_e_00} becomes, after integrating in time from $0$ to $1$:
\begin{equation}\label{t_b_1}
 \int_\Omega \tilde{g}^0_\ep(w_\ep,C_\ep)(1)  -  \int_\Omega \tilde{g}^0_\ep(w_\ep,C_\ep)(0) + \frac{1}{\alpha_\ep} \int_0^1  \int_\Omega \abs{\dt w_\ep}^2 +\abs{\dt C_\ep}^2   
=   \int_0^1 \int_\Omega V(w_\ep,C_\ep,0)  \cdot Z_\ep. 
\end{equation}
Note also that for any time-independent vector field $X\in C^0(\Omega;\Rn{2})$, 
\begin{multline}\label{t_b_9}
 \int_0^1 \int_\Omega V(w_\ep,C_\ep,0)  \cdot X 
= \int_0^1 \int_\Omega V(w_\ep,0,0)  \cdot X 
- \int_0^1 \int_\Omega \dt((1-\abs{w_\ep}^2)C_\ep \cdot X \\
=\int_0^1 \int_\Omega V(w_\ep,0,0)  \cdot X 
- \int_\Omega (1-\abs{w_\ep(1)}^2)C_\ep(1) \cdot X + \int_\Omega (1-\abs{w_\ep(0)}^2)C_\ep(0) \cdot X.
\end{multline}

Proposition \ref{diverge_control} implies that inequalities
\begin{equation}
 \int_{\partial\Omega} \abs{C_\ep(t) \cdot \nu}^2 \le  \alpha_\ep^2 j_{ex}^2 \pnorm{\nab f_1}{\infty}^2 \int_0^1 \int_{\partial \Omega} (1-\abs{w_\ep}^2)^2 
\end{equation}
and
\begin{equation}
 \int_\Omega \abs{\diverge{C_\ep(t)}}^2 \le 2  \int_0^1\int_\Omega \abs{\dt w_\ep}^2 + 2\alpha_\ep^2 \pnorm{f}{\infty}^2 \int_0^1 \int_\Omega (1-\abs{w_\ep}^2)^2
\end{equation}
both hold for all $t\in[0,1]$.  Condition \eqref{cond1'} guarantees that the hypotheses of Lemma \ref{bndry_conv} are satisfied by $w_\ep$ with a uniform constant $K$ for all time $t\in[0,1]$, so 
\begin{equation}
 \int_{\partial \Omega} (1-\abs{w_\ep}^2)^2 \le C \sqrt{\ep} \ale
\end{equation}
for a constant $C>0$ independent of $t$.  Hence the inequalities
\begin{equation}\label{t_b_4}
 \int_{\partial\Omega} \abs{C_\ep(t) \cdot \nu}^2 \le  C \alpha_\ep^2 j_{ex}^2 \pnorm{\nab f_1}{\infty}^2 \sqrt{\ep} \ale  
\end{equation}
and 
\begin{equation}\label{t_b_5}
 \int_\Omega \abs{\diverge{C_\ep(t)}}^2 \le 2  \alpha_\ep k_{ex} + C\alpha_\ep^2  \pnorm{f}{\infty}^2 \ep^2 \ale
\end{equation}
both hold for all $t\in[0,1]$.  We also clearly have that 
\begin{equation}\label{t_b_6}
 \int_\Omega \abs{\curl{C_\ep}(t)}^2 \le C \ale 
\end{equation}
for all $t\in[0,1]$.
We plug \eqref{t_b_4}--\eqref{t_b_6} into Proposition \ref{boundary_poincare} and recall that $f = j_{ex} f_1 - h_{ex} f_0$ to deduce the bound
\begin{equation}\label{t_b_7}
 \int_\Omega \abs{C_\ep(t)}^2 \le C(\ale + \alpha_\ep k_{ex}) + C \alpha_\ep^2 \sqrt{\ep} \ale k_{ex}^2 ( \pnorm{f_0}{\infty}^2 +  \pnorm{f_1}{\infty}^2+   \pnorm{\nab f_1}{\infty}^2).
\end{equation}

Note that 
\begin{equation}
 \alpha_\ep k_{ex} = \frac{\alpha_\ep}{\lep} k_{ex} \lep = o(1) \ale.
\end{equation}
When combined with the trivial bound $\abs{\nab w_\ep} \le \abs{\nab_{C_\ep}w_\ep} + \abs{C_\ep}$ and the bounds \eqref{t_b_7}, \eqref{cond1'}, this implies that
\begin{equation}\label{t_b_8}
 \hal \int_\Omega \abs{\nab w_\ep}^2 + \frac{(1-\abs{w_\ep}^2)^2}{2\ep^2}  \le C \ale
\end{equation}
for all $t\in[0,1]$.  This and condition \eqref{cond2'} then allow us to apply Proposition \ref{prod_est} to $w_\ep$ to deduce the convergence of $V(w_\ep,0,0)$ to some $V$ satisfying the bound \eqref{p_e_1}.  Plugging \eqref{t_b_8}, \eqref{cond2'}, and $\alpha_\ep k_{ex} \ll \ale$ into \eqref{p_e_1} then yields
\begin{equation}
 \abs{\int_0^1 \int_\Omega V \cdot X} \le 2 \frac{\pnorm{X}{\infty}}{\ale} \sqrt{C \ale \alpha_\ep k_{ex}} = o(1)
\end{equation}
for any $X \in C^0(\Omega \times [0,1])$, and hence that $V = 0$.  The vorticity measures $\mu(w_\ep,0)$  also converge to a limiting measure $\mu$, and the relation $\dt \mu + \curl{V}=0$ then implies that $\dt \mu = 0$.  We thus have that $\mu(t) = \mu(0)$ for all $t\in[0,1]$.  Moreover, the bounds  \eqref{cond1'} and \eqref{t_b_7} imply that $\mu(w_\ep,C_\ep)(t) \rightarrow \mu(t) = \mu(0)$ for each $t\in[0,1]$ as well, i.e. the vortices do not move.

Returning to \eqref{t_b_9}, we employ the convergence  $V(w_\ep,0,0)\rightarrow 0$  and the bounds \eqref{cond1'}, \eqref{t_b_7} to deduce that
\begin{equation}
 \int_0^1 \int_\Omega V(w_\ep,C_\ep,0) \cdot Z_\ep  = o(1) k_{ex}
\end{equation}
Plugging this into \eqref{t_b_1}, we see that
\begin{equation}\label{t_b_3}
 \int_\Omega \tilde{g}^0_\ep(w_\ep,C_\ep)(1)  -  \int_\Omega \tilde{g}^0_\ep(w_\ep,C_\ep)(0) \le o(1) k_{ex}
\end{equation}
We may therefore invoke  Lemma \ref{outside_bound} and Proposition \ref{lower_bound} to bound the modified energy at time $t=1$ from below.   Indeed, we find that
\begin{equation}\label{t_b_2}
\begin{split}
 \int_\Omega \tilde{g}_\ep(w_\ep,C_\ep)(1) & \ge  \pi n \ale + n \gamma +  \frac{1}{2} \int_\Omega \abs{Z_\ep}^2   +    W_{d}(a) + o(1)k_{ex}  \\
& \ge \int_\Omega \tilde{g}^0_\ep(w_\ep,C_\ep)(0) - C_0 k_{ex}  + o(1)k_{ex}.
\end{split}
\end{equation}
Inequality \eqref{t_b_3} implies that \eqref{cond1'} could not have been an equality at $t=1$.  Hence \eqref{cond2'} must have been an equality, i.e. 
 \begin{equation} 
 \frac{1}{\alpha_\ep} \int_0^{1} \int_\Omega \abs{\dt w_\ep}^2 + \abs{\dt C_\ep}^2   = 2C_0 k_{ex}.
\end{equation}
Plugging this and \eqref{t_b_2} back into \eqref{t_b_1} then yields the inequality
\begin{equation}
 -C_0 k_{ex} + o(1)k_{ex} + 2C_0 k_{ex} \le o(1) k_{ex},
\end{equation}
which is a contradiction as $\ep \rightarrow 0$ since $C_0>0$.  We deduce that it cannot be the case that $\liminf \alpha_\ep/\lep =0$, i.e. that there exists a constant $T_0>0$ so that $\alpha_\ep > T_0\lep$ as $\ep\rightarrow 0$.
\end{proof}

%%%%%%%%%%%%%%%%%%%%%%%%%%%%%%%%%%%%%%%%%%%%%%%%%%%%%%%%%%%%%%%%%%%%%%%%
\section{Limiting dynamics in the original time scale}\label{original_time_section}
%%%%%%%%%%%%%%%%%%%%%%%%%%%%%%%%%%%%%%%%%%%%%%%%%%%%%%%%%%%%%%%%%%%%%%%%

The following lemma allows us to remove singularities of the form $\{a_i\} \times [0,T]$ for solutions to parabolic equations with certain integral bounds.  It is used in the subsequent theorem for deriving the limiting equation for the phase excess.

\begin{lem}\label{removable_singularity}
 Suppose $u$ is a weak solution to the equation $\dt u - \Delta u = f$ on the set $B(0,r)\backslash \{0\} \times [0,T]$ for some $T,r>0$ and $f\in L^2(B(0,r)\times [0,T])$.  Further suppose that
\begin{equation}
 \sup_{0\le t \le T} \int_{B(0,r)} \abs{u(x,t)}^2 dx + \int_0^T \int_{B(0,r)} \abs{\nab_x u}^2 dx dt < \infty.
\end{equation}
Then $u$ may be extended to a solution in all of $B(0,r) \times [0,T]$.
\end{lem}
\begin{proof}
It is shown in \cite{alk} that under the above assumptions, $u$ may be extended to a solution of the equation on all of $B(0,r)\times [0,T]$ if and only if the thermal capacity of the set $\{0\} \times [0,T]$ vanishes.  That this set has vanishing thermal capacity is established in Corollary 1 of \cite{tayl_wat}.
\end{proof}

Define $\Theta_a$ to be the multi-valued function, harmonic on $\Omega \backslash \{a_i\}$, so that 
\begin{equation}\label{Theta_def}
 e^{i\Theta_a} = \prod_{i=1}^n \left(\frac{x-a_i}{\abs{x-a_i}} \right)^{d_i}.
\end{equation}
Note that while $\Theta_a$ is multi-valued, its gradient is well defined away from the points $\{a_i\}$.  We are now able to derive the limiting dynamics in the original time scale.  The following theorem, which is modeled on results in \cite{lin_1,spirn}, gives the dynamics for the pair $(v_\ep,B_\ep)$, and since $u_\ep = v_\ep e^{if}$ and $A_\ep = B_\ep + \nb h_0$, we may trivially derive the dynamics for $(u_\ep,A_\ep)$ from the theorem.

\begin{thm}\label{original_dynamics}
Suppose $j_{ex} =h_{ex}=1$ so that $J_{ex} = J$ and $H_{ex}=H$ do not depend on $\ep$.  Let $(v_\ep,B_\ep)$ solve \eqref{mod_eqn_1}--\eqref{mod_eqn_2} in the $\Phi=f$ gauge on $\Omega \times [0,T]$ for some $T>0.$  Further suppose that the initial data $(v_\ep(0),B_\ep(0))$ are well-prepared at order $C_0>0$ in the sense of \eqref{well_prepared_def}.  Then the following hold.
\begin{enumerate}
 \item\label{o_d_0_1} $\mu_\ep(t) \rightarrow \mu_0$ for all $t\in[0,T]$, i.e. the vortices do not move.
 
\item\label{o_d_0_2} $v_\ep \rightharpoonup v_*$ weakly in $H^1_{loc}(\Omega \backslash \{a_i\} \times [0,T])$, where 
\begin{equation}
 v_* = \prod_{i=1}^n \left(\frac{x-a_i}{\abs{x-a_i}}\right)^{d_i} e^{i\psi_*} = e^{i\Theta_a + i \psi_*},
\end{equation}
where $\Theta_a$ is defined by \eqref{Theta_def} and $\psi_*$ is a single-valued function on $\Omega\times [0,T]$ satisfying 
\begin{equation}
\begin{cases}
 \dt \psi_* - \Delta \psi_* + \psi_* = \psi_*(0) & \text{in } \Omega \\
 \nab \psi_* \cdot \nu = -\nab \Theta_a \cdot \nu & \text{on } \partial \Omega.
\end{cases}
\end{equation}

\item\label{o_d_0_3} $B_\ep \rightarrow B_*$  in $L^2(\Omega \times [0,T])$ and $h'_\ep \rightharpoonup h'_* = \curl{B_*}$ weakly in $L^2(\Omega \times [0,T])$.  The function $h'_*$ satisfies
\begin{equation}
\begin{cases}
 \dt h'_* - \Delta h'_* + h'_* = 2\pi \sum_{i=1}^n d_i \delta_{a_i} &  \text{in } \Omega \\
 h'_* = 0 & \text{on } \partial \Omega
\end{cases}
\end{equation}
with the PDE satisfied in $\mathcal{D}'(\Omega \times [0,T])$.

\end{enumerate}

\end{thm}

\begin{proof}
We divide the proof into several steps.\\
Step 1.

We begin by showing that the energy stays well behaved and that the vortices do not move in time. Since the initial data are well-prepared at order $C_0$ we may apply Theorem \ref{time_bound} to conclude that for $\ep$ sufficiently small (so that $T_0 \ale \ge T$) we have the bounds
\begin{equation}
 \int_\Omega \tilde{g}^0_\ep(v_\ep,B_\ep)(t) < \int_\Omega \tilde{g}^0_\ep(v_\ep,B_\ep)(0) + 2C_0 \text{ for all } t\in [0,T]
\end{equation}
and
\begin{equation}\label{o_d_003}
 \int_0^{T} \int_\Omega \abs{\dt v_\ep}^2 + \abs{\dt B_\ep}^2  < 2 C_0.
\end{equation}
These bounds imply, as in the proof of Theorem \ref{time_bound}, that the limiting velocity vanishes, $V=0$, and hence that $\mu_\ep(t) \rightarrow  \mu_0$ for all $t\in[0,T]$.  This proves item \ref{o_d_0_1}. \\\\
Step 2.

We now use the concentration of the energy around the vortex points $\{a_i\}$ to derive upper bounds for the energy away from the vortex points.  Fix 
\begin{equation*}
0 < \sigma < \frac{1}{2} \min\{ \dist(a_i,\partial \Omega)\} \cup \{\abs{a_i - a_j} \;\vert\; i\neq j  \}\cup\{6/\sqrt{2}\}
\end{equation*}
and define $\Omega_\sigma = \Omega \backslash \cup_{i=1}^n B(a_i,\sigma)$.  The bounds proved in Step 1 and the well-preparedness of the initial data allow us to apply Lemma \ref{outside_bound} to find that 
\begin{multline}\label{o_d_1}
 \hal \int_{\Omega_\sigma}\abs{\nab_{B_\ep} v_\ep}^2(t)  + \frac{1}{2\ep^2} (1-\abs{v_\ep}^2)^2(t) + \frac{1}{8} \int_\Omega \abs{\curl{B_\ep(t)}}^2 \\
\le  \pi n \log \frac{1}{\sigma} + n(\gamma +C) + W_{d}(a) + 3C_0.
\end{multline}
for every $t\in[0,T]$.

From \eqref{o_d_1} and Lemma \ref{diverge_control}, we see that $\curl{B_\ep}$ and $\diverge{B_\ep}$ are both bounded when measured in the $L^\infty([0,T];L^2(\Omega))$ norm. On the other hand, Lemma \ref{diverge_control} and Lemma \ref{bndry_conv} show that $B_\ep \cdot \nu \rightarrow 0$ in $L^\infty([0,T];L^2(\partial \Omega))$.  These facts, when combined with Proposition \ref{boundary_poincare} and the elliptic Hodge estimate
\begin{equation}
\norm{B}_{H^{1}(\Omega)} \le C\left(   \norm{\diverge{B}}_{L^2(\Omega)} + \norm{\curl{B}}_{L^2(\Omega)} + \norm{B\cdot \nu}_{H^{-1/2}(\partial\Omega)}    \right)
\end{equation}
show that
\begin{equation}\label{o_d_3}
 \sup_{0\le t \le T} \norm{B_\ep}_{H^{1}(\Omega)}  \le C_*
\end{equation}
for $C_*$ a constant that does not depend on $\ep$.  On the other hand, \eqref{o_d_3} implies that $\dt B_\ep$ is uniformly bounded in $L^2([0,T];L^2(\Omega))$.  This allows us to apply a result from \cite{simon} to deduce that, up to extraction, $B_\ep \rightarrow B_*$ in $L^2([0,T];H^{1/2}(\Omega))$.  In particular, we also have that $B_\ep \rightarrow B_*$ in $L^2(\Omega \times [0,T])$.

Since 
\begin{equation}
\abs{\nab v_\ep} \le  \abs{\nab_{B_\ep} v_\ep} + \abs{v_\ep B} 
\end{equation}
we also deduce that up to extraction $v_\ep \rightharpoonup v_*$ weakly in $H^1(\Omega_\sigma \times [0,T])$ and $v_\ep(t) \rightharpoonup v_*(t)$ weakly in $H^1(\Omega_\sigma)$ for each $t\in[0,T]$.  Clearly $v_*$ is unit valued.\\\\
Step 3.

%lifting in sobolev spaces reference?
We now further manipulate the bounds of the last step to derive the structure of $v_*$.  Since $v_* \in H^1(\Omega_\sigma)$, we may write $v_* = e^{i \varphi_*}$ for $\varphi_*$ a multi-valued function such that $\nab \varphi_*$ is well-defined and $\nab \varphi_* \in L^2(\Omega_\sigma)$.  Passing to the limit in \eqref{o_d_1} shows that 
\begin{equation}
 \hal \int_{\Omega_\sigma} \abs{\nab \varphi_*(t) - B_*(t) }^2 \le \pi n \log \frac{1}{\sigma} + n(\gamma +C) + W_{d}(a) + 3C_0 
\end{equation}
for every $t\in[0,T]$.  Define $Y_* = \nab \varphi_* - \nab \Theta_a,$ where $\Theta_a$ is defined by \eqref{Theta_def}.  Then a simple modification of standard arguments (cf. \cite{lin_1,spirn}) shows that 
\begin{equation}
 \hal \int_{\Omega_\sigma} \abs{\nab \varphi_*(t) - B_*(t) }^2  
= \hal \int_{\Omega_\sigma} \abs{Y_*(t) - B_*(t) }^2  + \pi n \log \frac{1}{\sigma} - O_\sigma(1).
\end{equation}
Hence
\begin{equation}
 \hal \int_{\Omega_\sigma} \abs{Y_*(t) - B_*(t)   }^2 \le n(\gamma +C)  + W_{d }(a) +  3C_0 +  O_\sigma(1).
\end{equation}
Since $B_*(t)$ is bounded in $L^2(\Omega)$, we may let $\sigma \rightarrow 0$ in this inequality to conclude that $Y_*(t)$ is well-defined and bounded in $L^2(\Omega)$.  

The convergence $\mu_\ep \rightarrow 2\pi \sum d_i \delta_{a_i}$, along with the convergence of $v_\ep$ and $B_\ep$ imply that 
\begin{equation}
 \curl{\nab \varphi_*} = 2\pi \sum_{i=1}^n d_i \delta_{a_i}
\end{equation}
in the sense of distributions.  From this and the definition of $\Theta_a$ we see that $\curl{Y_*(t)} = 0$ in the sense of distributions for each $t\in[0,T]$.  By the weak form of the Poincar\'{e} lemma (cf. Lemma 3 of \cite{bbm}), it then holds that $Y_* = \nab \psi_*$ for some $\psi_* \in H^1(\Omega)$.  Adjusting by a constant if necessary, for every $t\in[0,T]$ we may write $\psi_*(t) = \varphi_*(t) - \Theta_a$ in $\Omega_\sigma$.  Moreover, $\dt \psi_* = \dt \varphi_*$ in $\Omega_\sigma \times [0,T]$. \\\\
Step 4.

We now derive the equation satisfied by $\psi_*$.  Take $(iv_\ep,\cdot)$ with \eqref{mod_eqn_1}--\eqref{mod_eqn_2} and expand:
\begin{equation}\label{londonheat_precursor}
 (iv_\ep,\dt v_\ep) + \abs{v_\ep}^2 f = \diverge\left((iv_\ep,\nab v_\ep) -\abs{v_\ep}^2 B_\ep +  \abs{v_\ep}^2 Z_\ep  \right).
\end{equation}  
For each $\sigma>0$ the left side converges in $\mathcal{D}'(\Omega_\sigma \times [0,T])$ to $\dt \psi_* + f$, and the right side to $\Delta \psi_* - \diverge{B_*} + \Delta f$.  Then, by \eqref{f_def} $\Delta f =f$, so we have that $\dt \psi_* - \Delta \psi_* = -\diverge{B_*}$.  On the other hand, we know that $\dt \diverge{B_\ep} = (iv_\ep,\dt v_\ep) + (\abs{v_\ep}^2-1)f$, so we may pass to the limit to find that $\dt \diverge{B_*} = \dt \psi_*$.  Since $\diverge{B_*(0)}=0$, we find that $\diverge{A_*} = \psi_* - \psi_*(0)$.  We deduce that  
\begin{equation}
 \dt \psi_* - \Delta \psi_* + \psi_* = \psi_*(0)
\end{equation}
in $\Omega\backslash\{a_i\} \times[0,T]$.

It remains to show that the possible singularities at $\{a_i\}$ are actually removable.  To this end, we define the function $\eta_* = \psi_* - \bar{\psi}_*$, where 
\begin{equation*}
 \bar{\psi}_*(t) = \frac{1}{\abs{\Omega}} \int_\Omega \psi_*(t).
\end{equation*}
Then $\eta_*$ solves the equation $\dt \eta_* - \Delta \eta_* = \psi_*(0) - \psi_* - \dt \bar{\psi}_* \in L^2(\Omega \times [0,T]).$  By the Poincar\'{e} inequality
\begin{equation}\label{o_d_004}
 \int_\Omega \abs{\psi_*(t) - \bar{\psi}_*(t)}^2 \le C \int_\Omega \abs{\nab \psi_*(t)}^2 \le C_*
\end{equation}
for a constant $C_* = C(n,\Omega,h_{ex},d,a,C_0)$.  Hence
\begin{equation}
 \sup_{0\le t \le T} \int_\Omega \abs{\eta_*}^2 + \int_0^T \int_\Omega \abs{\nab \eta_*}^2 \le C_*(1+T).
\end{equation}
We may therefore apply Lemma \ref{removable_singularity} to conclude that the singularities at $\{a_i\}$ are removable and $\eta_*$ solves $\dt \eta_* - \Delta \eta_* = \psi_*(0) - \psi_* - \dt \bar{\psi}_*$ in all of $\Omega \times [0,T]$.  This then implies that $\psi_*$ solves \eqref{o_d_004} in all of $\Omega \times [0,T]$. The boundary condition $\nab \varphi_* \cdot \nu = 0$ carries over from the condition $\nab v_\ep \cdot \nu=0$ by multiplying \eqref{londonheat_precursor} by a test function that does not vanish on $\partial \Omega$ but that vanishes in a neighborhood of the vortex locations and passing to the limit.  So, $\nab \psi_* \cdot \nu = -\nab \Theta_a \cdot \nu$ on $\partial \Omega$. 
\\\\
Step 5.

We now derive the equation for $h'_* = \curl{B_*}$.  The magnetic potential $B_\ep$ satisfies the equation 
\begin{equation}\label{o_d_2}
\dt B_\ep  -\nab^\bot \curl{B_\ep} + B_\ep =  (iv_\ep,\nab_{B_\ep} v_\ep)+B_\ep + (\abs{v_\ep}^2-1) Z_\ep. 
\end{equation}
We may take the curl of \eqref{o_d_2} in the sense of distributions and pass to the limit, employing the convergence of $\mu_\ep$, to find that 
\begin{equation}
 \dt h'_* - \Delta h'_* + h'_*= 2\pi \sum_{i=1}^n d_i \delta_{a_i}
\end{equation}
in $\mathcal{D}'(\Omega \times [0,T])$.  That $h'_* = 0$ on $\partial \Omega$ carries over from the boundary condition $h'_\ep = \curl{B_\ep}=0$.

\end{proof}

%%%%%%%%%%%%%%%%%%%%%%%%%%%%%%%%%%%%%%%%%%%%%%%%%%%%%%%%%%%%%%%%%%%%%%%%
\section{Limiting dynamics in the accelerated time scale}\label{accelerated_time_section}
%%%%%%%%%%%%%%%%%%%%%%%%%%%%%%%%%%%%%%%%%%%%%%%%%%%%%%%%%%%%%%%%%%%%%%%%

%%%%%%%%%%%%%%%%%%%%%%%%%%%%%%%%%%%%%%%%%%%%%%%%%%%%%%%%%%%%%%%%%%%%%%%%
\subsection{Preliminaries and vortex motion}
%%%%%%%%%%%%%%%%%%%%%%%%%%%%%%%%%%%%%%%%%%%%%%%%%%%%%%%%%%%%%%%%%%%%%%%%

In this section we derive the limiting dynamics in the accelerated time scale.  We rescale in time at the scale $\lep:=\ale/k_{ex}$ with $k_{ex} = \max\{h_{ex},j_{ex}\}$ by making the substitutions
\eqref{acc_def}.  In this scaling the equations in the $\Phi=f$ gauge become
\begin{equation}\label{accelerated}
 \begin{cases}
   \lep^{-1} \dt v_\ep  = \Delta_{B_\ep} v_\ep + \frac{v_\ep}{\varepsilon^2}(1-\abs{v_\ep}^2) + 2i\nab_{B_\ep} v_\ep \cdot Z_\ep    - v_\ep  \abs{Z_\ep }^2 \\
  \lep^{-1}  \dt B_\ep  = \nab^\bot h'_\ep + (iv_\ep,\nab_{B_\ep} v_\ep) +(\abs{v_\ep}^2-1) Z_\ep
 \end{cases}
\end{equation}
along with the usual boundary and initial conditions.   In the accelerated scale, the evolution equation for the modified free energy density becomes
\begin{equation}\label{accelerated_en_evolve}
   \dt \tilde{g}^0_\ep(v_\ep,B_\ep) = \diverge(\dt v_\ep, \nab_{B_\ep} v_\ep) + \curl(h'_\ep \dt B_\ep)  - \frac{1}{\lep}\abs{\dt v_\ep}^2 - \frac{1}{\lep}\abs{\dt B_\ep}^2 + V(v_\ep,B_\ep,0)\cdot Z_\ep.
\end{equation}

We will assume throughout this section that the initial data $(v_\ep(0),B_\ep(0))$ are well-prepared at order $C_0 $ in the sense of \eqref{well_prepared_def}.  Writing 
\begin{equation}
\mu(0) = \sum_{i=1}^n d_i(0) \delta_{a_i(0)}
\end{equation}
for the $t=0$ limiting vortex measure, we further assume  that the initial vortex locations satisfy 
\begin{equation}
\min \{\abs{a_i(0) - a_j(0)} \;\vert\; i\neq j\} \cup \{\dist(a_i(0),\partial \Omega)\} \ge \sigma_0 
\end{equation} 
for some $\sigma_0>0$. We continue to assume that $h_{ex}, j_{ex}$ fall into one of the four regimes \eqref{regimes}.  In any of the four cases  $\lep \rightarrow \infty$ as $\ep \rightarrow 0$.

Rescaling in time at scale $\lep$, Theorem \ref{time_bound} provides for the existence of a constant $T_0>0$ such that 
\begin{equation}\label{acc_cond1}
 \int_\Omega \tilde{g}^0_\ep(v_\ep,B_\ep)(t) < \int_\Omega \tilde{g}^0_\ep(v_\ep,B_\ep)(0) + 2C_0 k_{ex}  \text{ for all } t\in [0,T_0 ]
\end{equation}
and 
\begin{equation}\label{acc_cond2}
 \int_0^{T_0 } \int_\Omega \abs{\dt v_\ep}^2 + \abs{\dt B_\ep}^2   < 2 C_0 k_{ex} \lep = 2 C_0 \ale.
\end{equation}
Note that the bound $k_{ex}\ll \ale^{1/9}$ and  \eqref{acc_cond2} imply that 
$\dt B_\ep/\lep$ and $\dt v_\ep/\lep$ vanish in $L^2(\Omega \times [0,T_0])$, and up to the extraction of a subsequence we may assume that $\dt B_\ep(t)/\lep$ and $\dt v_\ep(t)/\lep$ vanish in $L^2(\Omega)$ for almost every $t\in[0,T_0]$.  Recall also that from Lemma \ref{gradient_bound}, the bound $\pnorm{\nab_{B_\ep} v_\ep}{\infty} \le C/\ep$ holds for all time.

As the first order of business we record a lemma that shows the convergence of the space-time Jacobian and makes sense of the vortex trajectories.  The result also establishes for any $0<\sigma_*<\sigma_0$ the existence of a time $T_* = T_*(\sigma_*) \in(0,T_0]$ so that the vortex trajectories stay a distance $\sigma_*$ away from each other and $\partial \Omega$ for all $t\in[0,T_*]$.  We assume that such a $\sigma_*$ is fixed throughout the section, and we work exclusively in the domain $\Omega \times [0,T_*]$.

\begin{lem}\label{vortex_path}
Fix $0<\sigma_* < \sigma_0$. Then there exists a $T_* = T_*(\sigma_*)$ with $T_*\in(0,T_0]$ so that
\begin{enumerate}
\item The space-time Jacobian $(\mu(v_\ep,B_\ep),V(v_\ep,B_\ep,0)) \rightarrow (\mu,V)$ in $(C^{0,\alpha}(\Omega \times [0,T_0]))^*$ for all $\alpha\in(0,1)$.
\item There exist functions $a_i\in H^1([0,T_*];\Omega)$, $i=1,\dotsc,n$, so that for all  $t\in[0,T_*]$
\begin{equation}
\min \{\abs{a_i(t) - a_j(t)} \;\vert\; i\neq j\} \cup \{\dist(a_i(t),\partial \Omega)\} \ge \sigma_* 
\end{equation}
and
\begin{equation}
\mu(t) = 2\pi \sum_{i=1}^n d_i(0) \delta_{a_i(t)}.
\end{equation}
\end{enumerate}
\end{lem}

\begin{proof}
The bounds \eqref{acc_cond1} and \eqref{acc_cond2} allow us to argue as in the proof of Theorem \ref{time_bound} to deduce the first item.  For the second item, Proposition III.2 of \cite{ss_gamma} proves the existence of the vortex paths $a_i\in H^1([0,T_*];\Omega)$, and the embedding $H^1 \hookrightarrow C^{0,1/2}$ allows us to find the $T_* = T_*(\sigma_*)$ so that the vortex paths stay separated.  
\end{proof}

\begin{remark}
Since  the degree of the $i^{th}$ vortex does not change for any time in $[0,T_*]$, we may consolidate notation and write only $d_i$ in place of $d_i(0)$.
\end{remark}

According to Lemma \ref{vortex_path}, the functions $a_i\in C^{0,1/2}([0,T_*];\Omega)$ for $i=1,\dotsc,n$.  So, for any $0< \sigma < \sigma_*/4$, we may apply Lemma \ref{outside_bound} to find that 
\begin{equation}\label{d_l_1}
\int_{\Omega \backslash \cup B(a_i(t),\sigma)} \tilde{g}^0_\ep(v_\ep,B_\ep) \le k_{ex} C(\sigma,n,\Omega)
\end{equation}
and 
\begin{equation}
\int_\Omega \abs{h'_\ep(t)}^2 \le k_{ex} C(\sigma,n,\Omega)
\end{equation}
for all $t\in[0,T_*]$.  The latter bound implies that up to extraction 
\begin{equation}
\frac{h'_\ep(t)}{\sqrt{k_{ex}}} \rightharpoonup h'_*(t) \text{ weakly in } L^2(\Omega),
\end{equation}
whereas the former implies that 
\begin{equation}
\frac{h'_\ep(t)}{\sqrt{k_{ex}}} \rightharpoonup h'_*(t) \text{ weakly in }H^1(\Omega \backslash \cup B(a_i(t),\sigma)). 
\end{equation}
We may then integrate both sides of the equation
\begin{equation}
\frac{\dt B_\ep}{\lep} +(1-\abs{v_\ep}^2) Z_\ep = \nab^\bot h'_\ep + (iv_\ep,\nab_{B_\ep} v_\ep),
\end{equation}  
against $\nab^\bot \eta$ for $\eta \in C_c^\infty(\Omega)$, divide by $\sqrt{k_{ex}}$, integrate by parts, and pass to the limit to find that $h'_*(t)$ satisfies the PDE
\begin{equation}
-\Delta h'_*(t) + h'_*(t) = 
\begin{cases}
2\pi \sum_{i=1}^n d_i \delta_{a_i(t)} & \text{ if } k_{ex} =1 \\
0 & \text{ if } k_{ex} \gg 1 
\end{cases}
\end{equation}
in the sense of distributions.  By trace theory and the fact that $h'_\ep = 0$ on $\partial\Omega$, we also have that $h'_* =0$ on $\partial \Omega$.
Note that when $k_{ex}=1$, $h'_*$ is smooth outside of any neighborhood of $\cup_{i=1}^n a_i([0,T_*])$, and that when $k_{ex} \gg 1$, $h'_* =0$ identically.

%The energy upper bounds and Lemma \ref{outside_bound} show that $h_\ep$ is bounded in $L^2(\Omega \times %[0,T_*])$ so that up to extraction we have  $h_\ep \rightharpoonup h_*$ weakly in $L^2(\Omega \times %[0,T_*])$.
%We may then pass to the limit in the second equation of \eqref{accelerated}  to deduce that $h_*$ satisfies 
%\begin{equation}\label{h_star_def}
%-\Delta h_* + h_* =  2\pi \sum_{i=1}^n d_i(0) \delta_{a_i(t)}
%\end{equation}
%in the sense of distributions.

%%%%%%%%%%%%%%%%%%%%%%%%%%%%%%%%%%%%%%%%%%%%%%%%%%%%%%%%%%%%%%%%%%%%%%%%
\subsection{Convergence of the modified energy density}
%%%%%%%%%%%%%%%%%%%%%%%%%%%%%%%%%%%%%%%%%%%%%%%%%%%%%%%%%%%%%%%%%%%%%%%%

The term in the energy evolution equation that allows for the identification of the vortex locations at each time is the normalized energy density $\tilde{g}^0_\ep(v_\ep,B_\ep)(t)/\ale$, viewed as a measure on $\Omega$.  The mass of this measure is clearly bounded, so at any particular time we may extract a subsequence that converges in the weak sense of measures.  The technical obstruction is that this extracted subsequence depends on the choice of $t\in[0,T_*]$, whereas we would like the subsequence to converge for all $t\in[0,T_*]$.  Fortunately, the measures satisfy a certain semi-decreasing (in time) property that allows us to find such a subsequence by adapting results from \cite{bos_col}.  We begin with a proof of this semi-decreasing result.

\begin{lem}\label{semi_decreasing}
Let $\phi \in C^1(\Omega)$.  Then for any $0\le t_1 < t_2 \le T_*$ it holds that
 \begin{equation}\label{sem_d_0}
 \int_\Omega \phi^2 \frac{\tilde{g}^0_\ep(v_\ep,B_\ep)(t_2)}{\ale} - \int_\Omega \phi^2 \frac{\tilde{g}^0_\ep(v_\ep,B_\ep)(t_1)}{\ale} \le C(n,C_0,f,h_0,\phi) \sqrt{t_2-t_1} + o(1,\phi),
\end{equation}
where $C(n,C_0,f,h_0,\phi,\Omega)$ is a constant depending on $n, C_0, \pnorm{\abs{\nab f_1} + \abs{\nab f_0} + \abs{\nab h_0}}{\infty},$ and $\norm{\phi}_{C^1},$  and where $o(1,\phi)\rightarrow 0$ as $\ep \rightarrow 0$, with rate of convergence dependent on $\phi$ but independent of $t_1, t_2$.
\end{lem}
\begin{proof}
 Lemma \ref{mod_en_evolve}, rescaled to the accelerated time scale and then integrated in time from $t_1$ to $t_2$, shows that 
\begin{multline}\label{sem_d_3}
  \int_\Omega \phi^2 \frac{\tilde{g}^0_\ep(v_\ep,B_\ep)(t_2)}{\ale} -  \int_\Omega \phi^2 \frac{\tilde{g}^0_\ep(v_\ep,B_\ep)(t_1)}{\ale} 
+ \frac{1}{\lep \ale} \int_{t_1}^{t_2} \int_\Omega \phi^2 \left(\abs{\dt v_\ep}^2 + \abs{\dt B_\ep}^2 \right) \\
= I + II
\end{multline}
where
\begin{equation}
 I:= - \frac{1}{\ale}   \int_{t_1}^{t_2} \int_\Omega 2\phi \nab \phi \cdot (\dt v_\ep,\nab_{B_\ep} v_\ep) + 2\phi\nab^\bot\phi\cdot h'_\ep  \dt B_\ep 
\end{equation}
and 
\begin{equation}
II:=  \frac{1}{\ale} \int_{t_1}^{t_2} \int_\Omega \phi^2   V(v_\ep,B_\ep,0) \cdot Z_\ep. 
\end{equation}
We apply Cauchy-Schwarz and bounds \eqref{acc_cond1} and \eqref{acc_cond2} to $I$ to derive the inequality
\begin{multline}\label{sem_d_1}
\abs{I} \le \norm{\phi}_{C^1} \left( \int_{t_1}^{t_2} \int_\Omega \frac{\tilde{g}^0_\ep(v_\ep,B_\ep)}{\ale} \right)^{1/2}  \left( \frac{1}{\ale}\int_{t_1}^{t_2} \int_\Omega \abs{\dt v_\ep}^2 + \abs{\dt B_\ep}^2\right)^{1/2} \\ 
\le  C(n,C_0,\phi) \sqrt{t_2 - t_1}.
\end{multline}
To handle $II$ first note that, according to Lemma \ref{vortex_path}, for a fixed vector field $X\in C^0(\Omega;\Rn{2})$
\begin{equation}
 \int_{t_1}^{t_2} \int_\Omega \phi^2 X \cdot V(v_\ep,B_\ep,0) = \int_{t_1}^{t_2} \int_\Omega \phi^2 V \cdot X + o(1,\phi),
\end{equation}
where $o(1,\phi)\rightarrow 0$ as $\ep \rightarrow 0$, with rate of convergence dependent on $\phi$ but independent of $t_1, t_2$.  We then apply item $2$ of Proposition \ref{prod_est} along with the bounds \eqref{acc_cond1} and \eqref{acc_cond2} to bound
\begin{equation}
 \abs{ \int_{t_1}^{t_2} \int_\Omega \phi^2 X \cdot V } \le C(n,C_0,\phi,X).\sqrt{t_2 - t_1}.
\end{equation}
Applying this with $X=\nab f_1, \nab f_0, \nb h_0$ and employing the fact that $k_{ex} \ll \ale$ then shows that 
\begin{equation}\label{sem_d_2}
 \abs{II} \le o(1)\sqrt{t_2 - t_1} + o(1,\phi).
\end{equation}
Then \eqref{sem_d_0} follows by plugging the bounds \eqref{sem_d_1} and \eqref{sem_d_2} into \eqref{sem_d_3}.
\end{proof}

With this semi-decreasing property established, we can now show that up to the extraction of a \emph{single} subsequence, the modified energy density converges to a sum of Dirac masses \emph{for all} $t\in[0,T_*]$.

\begin{prop}\label{density_convergence}
There exists a subsequence so that
\begin{equation}\label{den_con_0}
 \frac{\tilde{g}^0_\ep(v_\ep,B_\ep)(t)}{\ale} \wstar \pi \sum_{i=1}^n \delta_{a_i(t)} 
\end{equation}
weakly-$*$ in $(C^1(\Omega))^*$ for all $t\in[0,T_*]$.

\end{prop}
\begin{proof}
We will first establish the convergence of the measures and then establish the structure of the limiting measure.  In order to establish the convergence (up to extraction) for all time, we will use a variant of Helly's selection principle in conjunction with the semi-decreasing property proved in Lemma \ref{semi_decreasing}.  We essentially follow the strategy presented in Section 5.4 of \cite{bos_col}.

Let $\{\phi_k\}_{k=1}^\infty$ be a countable set of functions in $C^1(\Omega)$ so that the span of $\{\phi_k^2\}_{k=1}^\infty$ is dense.  For each $k,\ep$ define the function $\xi_{k,\ep}:[0,T_*] \rightarrow \Rn{}$ by 
\begin{equation}
 \xi_{k,\ep}(t) = \int_\Omega \phi_k^2 \frac{\tilde{g}^0_\ep(v_\ep,B_\ep)(t)}{\ale}.
\end{equation}
By Lemma \ref{semi_decreasing}, the functions $\xi_{k,\ep}$ satisfy the following semi-decreasing property:
for every $\delta>0$ there exist $\ep_k>0$ and $\tau_k>0$ so that for every $t_2\in(0,T)$ and $t_1 \in (t_2 - \tau_k,t_2)$ it holds that
\begin{equation}
 \xi_{k,\ep}(t_2) \le \xi_{k,\ep}(t_1) + \delta \; \text{ for all } \ep< \ep_k.
\end{equation}
Then a semi-decreasing variant of Helly's selection theorem (cf. Lemma 5.4 in \cite{bos_col}) implies that there exists a set of functions $\xi_k:[0,T_*]\rightarrow \Rn{}$ such that up to extraction 
\begin{equation}
\xi_{k,\ep}(t) \rightarrow \xi_k(t) \;\text{ for all }t\in[0,T_*] \text{ and for all } k\in\mathbb{N}
\end{equation}
as $\ep \rightarrow 0$.  From the density of the span of $\{\phi_k^2\}_{k=1}^\infty$, we deduce that there exists a family of measures $\nu(t)$ such that 
\begin{equation}
 \frac{\tilde{g}^0_\ep(v_\ep,B_\ep)(t)} {\ale} \wstar \nu(t) \;\text{ weakly-}* \text{ in }(C^1(\Omega))^* \text{ for all } t\in[0,T_*].
\end{equation}

We now derive the structure of the limiting measures $\nu(t)$.  According to Lemmas \ref{vortex_path} and \ref{outside_bound}, for $\sigma < \sigma_*/4$ 
\begin{equation}
 \int_{\Omega \backslash \cup B(a_i(t),\sigma)} \frac{\tilde{g}^0_\ep(v_\ep,B_\ep)(t) }{\ale} \le o(1).
\end{equation}
Since $\sigma$ can be taken to be arbitrarily small, this then implies that 
\begin{equation}
 \nu(t) = \sum_{i=1}^n \alpha_i(t) \delta_{a_i(t)}.
\end{equation}
We now calculate the value of $\alpha_i(t)$.  Fixing $\sigma < \sigma_*/4$, for each $i$ we may choose $\eta_i \in C^1(\Omega)$ so that $\supp(\eta_i)=B(a_i(t),2\sigma)$, $\eta_i =1$ on $B(a_i(t),\sigma)$, and $0\le \eta_i \le 1$.  Then applying Lemma \ref{inverse_measure} with this choice of $\sigma$, we have that 
\begin{equation}
 \alpha_i(t) = \lim_{\ep\rightarrow 0} \int_\Omega \eta_i \frac{\tilde{g}^0_\ep(v_\ep,B_\ep)(t)}{\ale} 
 \ge \liminf_{\ep \rightarrow 0} \int_{B(a_i(t),\sigma)}  \frac{\tilde{g}^0_\ep(v_\ep,B_\ep)(t)}{\ale} \ge \pi.
\end{equation}
On the other hand, from the bound \eqref{acc_cond1} it holds that 
\begin{equation}
 \sum_{i=1}^n \alpha_i(t) = \lim_{\ep \rightarrow 0} \sum_{i=1}^n \int_\Omega \eta_i \frac{\tilde{g}^0_\ep(v_\ep,B_\ep)(t)}{\ale} \\
\le \limsup_{\ep\rightarrow 0} \int_\Omega \frac{\tilde{g}^0_\ep(v_\ep,B_\ep)(t)}{\ale} \le \pi n.
\end{equation}
Hence $\alpha_i(t) =\pi$ for each $i=1,\dotsc,n$.

\end{proof}

\begin{remark}
 We assume in what follows that we are working with the extracted subsequence so that the convergence result of Proposition \ref{density_convergence} holds.
\end{remark}

%%%%%%%%%%%%%%%%%%%%%%%%%%%%%%%%%%%%%%%%%%%%%%%%%%%%%%%%%%%%%%%%%%%%%%%%
\subsection{Convergence of the stress-energy tensor}
%%%%%%%%%%%%%%%%%%%%%%%%%%%%%%%%%%%%%%%%%%%%%%%%%%%%%%%%%%%%%%%%%%%%%%%%

The stress-energy tensor associated to a configuration $(v_\ep,B_\ep)$ is the symmetric 2-tensor, $T_\ep$, with components 
\begin{equation}\label{stress_def}
(T_\ep)_{ij} = (\partial_i^{B_\ep} v_\ep, \partial_j^{B_\ep} v_\ep) - \delta_{ij} \left(\hal \abs{\nab_{B_\ep} v_\ep}^2 + \frac{1}{4\varepsilon^2}(1-\abs{v_\ep}^2)^2 - \hal (\curl{B_\ep})^2 \right),
\end{equation}
where $\partial_i^{B} v := \partial_i v - i B_i v$ is the $i^{th}$ covariant partial derivative.  The divergence of $T_\ep$ is the vector $\diverge{T_\ep}$ with components 
\begin{equation}\label{divergence_def}
(\diverge{T_\ep})_i = \partial_1 T_{1i} + \partial_2 T_{2i}.
\end{equation}
The divergence encodes the ``force'' acting on the vortices, and by passing to the limit in $T_\ep$ for solutions $(v_\ep,B_\ep)$ we will be able to derive one of the terms driving the limiting vortex dynamics.  The method of passing to the limit in the stress-energy tensor has been used extensively in the study of Ginzburg-Landau dynamics \cite{js_2,lin_1,lin_2,spirn,serf_2,bos_col}.

To make sense of the limit of $T_\ep$, we  follow a strategy similar to that of Chapter 13 of \cite{ss_book}, where they study the limit of $T_\ep$ for $(v_\ep,B_\ep)$ solutions to the elliptic Ginzburg-Landau equations on $\Omega$ (critical points of the Ginzburg-Landau energy functional).  Roughly, the idea is to use the vanishing of the right hand side of the equation
\begin{equation*}
\nab^\bot h'_\ep + (iv_\ep,\nab_{B_\ep} v_\ep) = \frac{\dt B_\ep}{\lep} + (1-\abs{v_\ep}^2)Z_\ep
\end{equation*}
to show that $T_\ep$ has the same limit as a similar tensor defined in terms of the induced magnetic field $h_\ep$, and then to derive the structure of the limit of the latter tensor.  In \cite{ss_book}, they prove a convergence in finite-parts result for $T_\ep$ on all of $\Omega$, which is stronger than what we shall prove here for $\Omega \times [0,T_*]$.  Indeed, the convergence result we prove here  holds only for cylinders $U_r \times [t_1,t_2]$ for which we know certain strong bounds on the energy.  Here we have written $U_r$ for a ball of radius $r$ rather than $B_r$ to avoid confusion with the vector potential $B_\ep$.

\begin{thm}\label{tensor_converge}
Suppose that in a ball of radius $r>0$, $U_r \subset \Omega$, it holds that
\begin{equation}\label{ten_con_0}
 \int_{U_r} \tilde{g}_\ep(v_\ep,B_\ep)(t) \le C k_{ex} \text{ for all } t\in[t_1,t_2].
\end{equation}
For the limiting induced magnetic field, $h'_*$, define the tensor 
\begin{equation}
 S(h'_*) := -\nab h'_* \otimes \nab h'_* + I_{2\times2} \left(\hal \abs{\nab h'_*}^2 + \hal \abs{h'_*}^2 \right).
\end{equation}
Then for $0 < s_1 < r$, 
\begin{equation}
  \frac{1}{k_{ex}} T_\ep \rightarrow S(h_*) \text{ in } L^1(U_{s_1} \times [t_1,t_2]).
\end{equation}
\end{thm}

\begin{remark}
 Recall that if $k_{ex} \gg 1$, then the limiting magnetic field vanishes, i.e. $h'_*=0$.  So, in this case, the theorem says that $T_\ep / k_{ex} \rightarrow 0$ in $L^1(U_{s_1} \times [t_1,t_2])$.
\end{remark}

The proof of the theorem is based on the following three lemmas, all of which rely heavily on the 
energy bound \eqref{ten_con_0}.  We begin by showing that $\abs{v_\ep}$ must be close to $1$.  The argument is a modification of one used in \cite{lin_2, spirn}.

\begin{lem}\label{modulus_close}
Fix $0<s_1<s_2<r$.  Then for almost every $t\in[t_1,t_2]$ it holds that
\begin{equation}
\lim_{\ep \rightarrow 0} \lep^{1/8} \pnormspace{1-\abs{v_\ep(t)}}{\infty}{U_{s_2}} =0
\end{equation}
\end{lem}
\begin{proof}
Fix $t\in[t_1,t_2]$ so that $\dt v_\ep(t)/\lep \rightarrow 0$, $\dt B_\ep(t) /\lep \rightarrow 0$ in $L^2(\Omega)$.  We know that the set of such $t$ has full measure.  For the rest of the proof we will neglect to write the dependence of the functions on time, but all are implicitly evaluated at the chosen time $t$. We assume that $\ep < (r - s_2)^2/k_{ex}^2$ so that $B(x,k_{ex} \sqrt{\ep}) \subset U_r$ for all $x\in U_{s_2}$.  Suppose, by way of contradiction, that 
\begin{equation}
\limsup_{\ep \rightarrow 0}  \lep^{1/8} \pnormspace{1-\abs{v_\ep}}{\infty}{U_{s_2}}  >0. 
\end{equation}
We may extract a subsequence (still denoted by $\ep$) so that 
\begin{equation}
 \lep^{1/8} \pnormspace{1-\abs{v_\ep}}{\infty}{U_{s_2}} \ge \alpha >0 \text{ for all } \ep.
\end{equation}
Since each $v_\ep$ is continuous, we may choose points $x_\ep  \in U_{s_2}$ with the property that  $\abs{v_\ep(x_\ep)} \le 1 - \alpha \lep^{-1/8}$.

Rewrite the equations \eqref{accelerated} in elliptic form as
\begin{equation}\label{mod_c_2}
\begin{cases}
\Delta_{B_\ep} v_\ep + \frac{v_\ep}{\ep^2}(1-\abs{v_\ep}^2) = K_\ep \\
\nab^\bot h'_\ep + (iv_\ep,\nab_{B_\ep}v_\ep) = X_\ep
\end{cases}
\end{equation}
where
\begin{equation}
K_\ep := \frac{\dt v_\ep}{\lep} -2 i\nab_{B_\ep} v_\ep \cdot Z_\ep + v_\ep  \abs{Z_\ep}^2 \text{ and } X_\ep := \frac{\dt B_\ep}{\lep} +(1-\abs{v_\ep}^2) Z_\ep. 
\end{equation}
Since the quantities of interest are gauge invariant, we are free to switch to the Coulomb gauge (again, only at the time $t$) via $v_\ep \mapsto w_\ep:= v_\ep e^{i\xi_\ep}$, $B_\ep \mapsto C_\ep:= B_\ep + \nab \xi_\ep$ with $\xi_\ep$ chosen so that
\begin{equation}
\begin{cases}
\diverge{C_\ep} = 0 & \text{in } U_{r} \\
C_\ep \cdot \nu = 0 & \text{on } \partial U_r.
\end{cases}
\end{equation}
Then the elliptic equations satisfied by $(w_\ep,C_\ep)$ are
\begin{equation}\label{mod_c_3}
\begin{cases}
\Delta_{C_\ep} w_\ep + \frac{w_\ep}{\ep^2}(1-\abs{w_\ep}^2) = e^{i \xi_\ep}K_\ep \\
\nab^\bot h'_\ep + (iw_\ep,\nab_{C_\ep}w_\ep) = X_\ep
\end{cases}
\end{equation}
where $k_\ep$, $X_\ep$ are as above (still with $v_\ep$ in their definition) and $h'_\ep = \curl{B_\ep} = \curl{C_\ep}$.  In the Coulomb gauge, we have that (cf. Proposition 3.3 of \cite{ss_book}) 
\begin{equation}
\norm{C_\ep}^2_{H^2(U_r)} \le C \norm{h'_\ep}^2_{H^1(U_r)} 
\le C\left(\int_{U_r} \tilde{g}^0_\ep(v_\ep,B_\ep) + \abs{\frac{\dt B_\ep}{\lep}}^2 + (1-\abs{v_\ep}^2)^2\abs{Z_\ep}^2 \right),
\end{equation}
which implies that $\pnormspace{C_\ep}{\infty}{U_r} \le C k_{ex}$.

Note that
\begin{equation}
 \int_{B_r} \tilde{g}^0_\ep(w_\ep,C_\ep) \ge \int_{B(x_\ep,k_{ex} \sqrt{\ep})} \tilde{g}^0_\ep(w_\ep,C_\ep) =  \int_{k_{ex} \ep}^{k_{ex} \sqrt{\ep}}  \left(s\int_{\partial B(x_\ep,s)} \tilde{g}^0_\ep(w_\ep,C_\ep)\right) \frac{ds}{s},
\end{equation}
which implies that
\begin{equation}
\inf_{k_{ex} \ep < s < k_{ex} \sqrt{\ep}} \left(s \int_{\partial B(x_\ep,s)} \tilde{g}^0_\ep(w_\ep,C_\ep) \right)\le \frac{2}{\ale} \int_{U_r} \tilde{g}^0_\ep(w_\ep,C_\ep) \le \frac{2C k_{ex}}{\ale}=\frac{2C}{\lep},
\end{equation}
so that there exists $r_\ep \in (k_{ex} \ep, k_{ex} \sqrt{\ep})$ with the property that
\begin{equation}\label{mod_c_1}
r_\ep \int_{\partial B(x_\ep,r_\ep)} \tilde{g}^0_\ep(w_\ep,C_\ep)  \le \frac{4C}{\lep}.
\end{equation}
We now recall the Pohozaev identity on $B(x_\ep,r_\ep)$ for solutions to \eqref{mod_c_3}:
\begin{multline}
\int_{B(x_\ep,r_\ep)} \frac{1}{2\varepsilon^2}(1-\abs{w_\ep}^2)^2 - (h'_\ep)^2 
=  \int_{B(x_\ep,r_\ep)} ( (e^{i\xi_\ep} K_\ep,\nab_{C_\ep }w_\ep) - h'_\ep X_\ep^\bot )\cdot (x-x_\ep) \\
+ r_\ep \int_{\partial B(x_\ep,r_\ep)} \hal \abs{\nab_{C_\ep} w_\ep \cdot \tau}^2 - \hal \abs{\nab_{C_\ep} w_\ep \cdot \nu}^2
  + \frac{1}{4\varepsilon^2}(1-\abs{w_\ep}^2)^2 -   \hal (h'_\ep)^2.
\end{multline}
We can use the $L^2(U_r)$ bounds on $K_\ep$ and $X_\ep$ that come from \eqref{ten_con_0} to control the first term on the right side, and \eqref{mod_c_1} controls the second term:
\begin{equation}
\int_{B(x_\ep,r_\ep)} \frac{1}{2\varepsilon^2}(1-\abs{w_\ep}^2)^2 - (h'_\ep)^2 
\le   C k_{ex}^3 r_\ep  + \frac{4C}{\lep} \le C k_{ex}^4 \sqrt{\ep} +\frac{4C}{\lep}.
\end{equation}
Now, by Sobolev 
\begin{equation}
\int_{B(x_\ep,r_\ep)} (h'_\ep)^2 \le C r_\ep \pnormspace{h'_\ep}{4}{U_r}^2 \le C k_{ex} \sqrt{\ep} \norm{h'_\ep}^2_{H^1(U_r)} \le C k_{ex}^4 \sqrt{\ep}.
\end{equation}
Hence 
\begin{equation}\label{mod_c_4}
\int_{B(x_\ep,r_\ep)} \frac{1}{2\varepsilon^2}(1-\abs{w_\ep}^2)^2 \le  C k_{ex}^4 \sqrt{\ep} +\frac{4C}{\lep} \le \frac{C}{\lep}
\end{equation}
for $\ep$ sufficiently small.

On the other hand, we have that 
\begin{equation}
\abs{\nab w_\ep} \le \abs{\nab_{C_\ep} w_\ep} +  \abs{C_\ep} = \abs{\nab_{B_\ep} v_\ep} +  \abs{C_\ep} \le  \frac{C}{\ep} +C k_{ex} \le \frac{C_0}{\ep} .
\end{equation}
Let $s_\ep =  \alpha(2C_0)^{-1} \ep \lep^{-1/8}$, and note that $s_\ep < r_\ep$ for $\ep$ sufficiently small.  Then the mean value theorem shows that 
\begin{equation}
\frac{\alpha}{2\lep^{1/8}} \le 1-\abs{w_\ep(x)}^2 \text{ for all } x\in B(x_\ep,s_\ep),
\end{equation}
and hence
\begin{equation}\label{mod_c_5}
\frac{\pi \alpha^4}{32 C_0^2 \lep^{1/2} } \le \int_{B(x_\ep,s_\ep)} \frac{(1-\abs{w_\ep}^2)^2}{2\ep^2}.
\end{equation}
Comparing \eqref{mod_c_4} and \eqref{mod_c_5}, we deduce that
\begin{equation}
\frac{\pi \alpha^4 }{32 C_0^2 \lep^{1/2} } \le \frac{C}{\lep},
\end{equation}
which yields a contradiction as $\ep \rightarrow 0$.  

\end{proof}

The next lemma, which is a parabolic modification of the elliptic result in Proposition 13.4 of \cite{ss_book}, shows that we can essentially replace $\nab_{B_\ep} v_\ep$ with $(iv_\ep,\nab_{B_\ep} v_\ep)$ in the definition of $T_\ep$.  Note that this lemma is the only place we use the full strength of the upper bound $k_{ex} \ll \ale^{1/9}$.

\begin{lem}\label{tensor_swap}
Define the tensor
\begin{equation}
 T_\ep' = (iv_\ep,\nab_{B_\ep} v_\ep) \otimes (iv_\ep,\nab_{B_\ep} v_\ep) - I_{2\times2} \left(\hal \abs{(iv_\ep,\nab_{B_\ep} v_\ep)}^2 -\hal (h'_\ep)^2\right). 
\end{equation}
Then 
\begin{equation}\label{t_swap_0}
 k_{ex}^{-1}\pnormspace{T_\ep - T_\ep'}{1}{U_{s_1} \times [t_1,t_2]} \rightarrow 0.
\end{equation}
\end{lem}
\begin{proof}
By Lemma \ref{modulus_close}, for almost every $t\in[t_1,t_2]$, as $\ep \rightarrow 0$, we may write $v_\ep(t) = \rho_\ep(t) e^{i\varphi_\ep(t)}$ for $\rho_\ep$, $\varphi_\ep$ well-defined and single-valued in $U_{s_2}$. We again neglect to write the dependence on $t$ in what follows.  Then since $(iv_\ep,\nab_{B_\ep} v_\ep) = \rho_\ep^2(\nab \varphi_\ep - B_\ep)$, 
\begin{equation}
\rho_\ep^2 T_\ep - T_\ep' = \rho_\ep^2\left(\nab \rho_\ep \otimes \nab \rho_\ep - I_{2\times2} \left(\hal \abs{\nab \rho_\ep}^2 + \frac{(1-\rho_\ep^2)^2}{4\ep^2} \right) \right) + I_{2\times2} (\rho_\ep^2 -1)\frac{(h'_\ep)^2}{2}.
\end{equation}
We may then write
\begin{equation}
 T_\ep - T_\ep' = (1-\rho_\ep^2)T_\ep + \rho_\ep^2 T_\ep - T_\ep'
\end{equation}
so that
\begin{equation}
 \int_{U_{s_1}}\frac{\abs{T_\ep - T_\ep'}}{k_{ex}} \le \frac{C}{k_{ex}}\int_{U_{s_1}}\abs{1-\rho_\ep^2} \tilde{g}^0_\ep(v_\ep,B_\ep)  
+\frac{C}{k_{ex}}\int_{U_{s_1}}\hal \abs{\nab \rho_\ep}^2 + \frac{(1-\rho_\ep^2)^2}{4\ep^2}.
\end{equation}
The first term on the right is easy to manage in view of \eqref{ten_con_0} and Lemma \ref{modulus_close}:
\begin{equation}
\frac{C}{k_{ex}}\int_{U_{s_1}}\abs{1-\rho_\ep^2} \tilde{g}^0_\ep(v_\ep,B_\ep)  \rightarrow 0.
\end{equation}
To handle the second term, we take $(v_\ep,\cdot)$ with the first equation of \eqref{accelerated} to find that  $\rho_\ep$ satisfies the PDE
\begin{equation}\label{t_swap_1}
 \frac{\dt \rho_\ep}{\lep} - \Delta \rho_\ep + \rho_\ep \abs{\nab \varphi_\ep +  Z_\ep - B_\ep}^2 = \frac{\rho_\ep}{\ep^2}(1-\rho_\ep^2).
\end{equation}
Multiply \eqref{t_swap_1} by $(1-\rho_\ep)$ and integrate over a ball $U_s$ for $0<s<s_2$ to be chosen later.  After integrating by parts and rearranging, we arrive at the equation
\begin{multline}\label{t_swap_2}
 \int_{U_s} \abs{\nab \rho_\ep}^2 + \frac{\rho_\ep}{\ep^2}(1-\rho_\ep^2)(1-\rho_\ep) = \int_{U_s} \rho_\ep (1-\rho_\ep) \abs{\nab \varphi_\ep +   Z_\ep - B_\ep}^2 \\ 
+ \int_{U_s}\frac{\dt \rho_\ep}{\lep} (1-\rho_\ep) - \int_{\partial U_s} \frac{\partial \rho_\ep}{\partial \nu} (1-\rho_\ep)
\end{multline}
By a mean-value argument, we may choose $s\in[s_1,s_2]$ so that 
\begin{equation}
 \int_{\partial U_s} \abs{\nab \rho_\ep}^2 \le \frac{2}{s_2-s_1} \int_{U_r} \tilde{g}^0_\ep(v_\ep,B_\ep)\le  \frac{C k_{ex}}{s_2-s_1}.
\end{equation}
On $U_{s_2},$ for $\ep$ sufficiently small,  it holds that
\begin{equation}\label{t_swap_3}
 \frac{1}{4} \le \frac{\rho_\ep}{1+\rho_\ep} \;\;\text{ and }\;\; \frac{1-\rho_\ep}{\rho_\ep} \le \frac{C}{\lep^{1/8}}.
\end{equation}
Combining \eqref{t_swap_2}--\eqref{t_swap_3}, we find that
\begin{multline}
 \frac{1}{k_{ex}} \int_{U_s} \hal \abs{\nab \rho_\ep}^2 + \frac{(1-\rho_\ep^2)^2}{4\ep^2} \le 
\frac{C}{k_{ex} \lep^{1/8}}  \int_{U_s} \abs{\nab_{B_\ep} v_\ep +  iv_\ep Z_\ep }^2 \\ 
+ \frac{C}{k_{ex} \lep^{1/8}} \left( \sqrt{\pi}s \left(\int_{U_s } \abs{\frac{\dt v_\ep}{\lep}}^2  \right)^{1/2}   + \sqrt{2\pi s} \sqrt{ \frac{C k_{ex} }{s_2-s_1}}  \right) 
\le \frac{C k_{ex}}{\lep^{1/8}} = o(1),
\end{multline}
where the last equality follows from the fact that $\lep = \ale/k_{ex}$ and $k_{ex}\ll \ale^{1/9}$.

Hence
\begin{equation}
\frac{1}{k_{ex}}\int_{U_{s_1}} \abs{T_\ep - T_\ep'}(t) \rightarrow 0 \text{ for almost every } t\in[t_1,t_2].
\end{equation}
Since for all $t\in[t_1,t_2]$,
\begin{equation}
\frac{1}{k_{ex}} \int_{U_{s_1}} \abs{T_\ep - T_\ep'}(t) \le \frac{C}{k_{ex}} \int_{U_r} \tilde{g}^0_\ep(v_\ep,B_\ep) \le C,
\end{equation}
we conclude that \eqref{t_swap_0} holds via an application of the dominated convergence theorem.

\end{proof}

The third lemma establishes the strong convergence of $h'_\ep/\sqrt{k_{ex}}$ and $\nab h'_\ep /\sqrt{k_{ex}}$ via an argument like that of Claim 5.4 in \cite{spirn}.

\begin{lem}\label{h_converge} It holds that
\begin{equation}
 \frac{h'_\ep}{\sqrt{k_{ex}}} \rightarrow h'_* \text{ and } \frac{\nab h'_\ep}{\sqrt{k_{ex}}} \rightarrow \nab h'_* \text{ in } L^2(U_{s_1} \times [t_1,t_2]).
\end{equation} 
\end{lem}
\begin{proof}
 Since $\nab^\bot h'_\ep = -(iv_\ep,\nab_{B_\ep}v_\ep) + \dt B_\ep/\lep + (1-\abs{v_\ep}^2) Z_\ep$, we know that 
\begin{equation}\label{h_con_1}
 \frac{1}{k_{ex}} \int_{U_r} \abs{\nab h'_\ep}^2 \le \frac{C}{k_{ex}} \int_{u_r} \abs{\nab_{B_\ep} v_\ep}^2 +  (1-\abs{v_\ep}^2)^2\abs{ Z_\ep }^2 + \abs{\frac{\dt B_\ep}{\lep}}^2,
\end{equation}
so for almost every $t\in[t_1,t_2]$ it holds that 
\begin{equation}
\frac{1}{k_{ex}} \int_{U_r} \abs{h'_\ep}^2 + \abs{\nab h'_\ep}^2(t) \le C.
\end{equation}
%this shows that there are subsequential limits, but these limits should all agree with the limit on the cylinder, $h_*$, so that the weak limit exists without subsequential extraction (limsup and liminf both are h_*)%
Then $h'_\ep/\sqrt{k_{ex}} \rightharpoonup h_*$ weakly in $H^1(U_r)$,  $h'_\ep/\sqrt{k_{ex}} \rightarrow h'_*$ in $L^2(U_r)$, and by dominated convergence $h'_\ep/\sqrt{k_{ex}} \rightarrow h_*$ in $L^2(U_s \times [t_1,t_2])$ as well.  Write $v_\ep = \rho_\ep e^{i\varphi_\ep}$ in $U_{s_2}$.  Then
\begin{equation}
 \frac{1}{\rho_\ep^2} \frac{\dt B_\ep}{\lep} + \frac{1-\rho_\ep^2}{\rho_\ep^2}  Z_\ep - \nab \varphi_\ep + B_\ep = \frac{1}{\rho_\ep^2} \nab^\bot h'_\ep .
\end{equation}
Fix $\xi \in C_c^\infty(U_{s_2})$.  Multiply the last equation by $\nab^\bot(\xi(h'_\ep-h'_*))$ and integrate over $U_r$ to get
\begin{multline}
 \int_{U_r} \frac{\xi}{\rho_\ep^2} \left(\abs{\nab h'_\ep}^2 - \nab h'_\ep \cdot \nab h'_* \right) = -\int_{U_r} \frac{h'_\ep-h'_*}{\rho_\ep^2}\nab^\bot \xi \cdot \nab^\bot h'_\ep + \xi h'_\ep (h'_\ep - h'_*) \\
+ \int_{U_r}\left(\frac{1}{\rho_\ep^2} \frac{\dt B_\ep}{\lep} + \frac{1-\rho_\ep^2}{\rho_\ep^2} Z_\ep \right)\cdot \nab^\bot(\xi(h'_\ep-h'_*)).
\end{multline}
From this and the known convergence results, we deduce that for almost every $t\in[t_1,t_2]$,
\begin{equation}
\frac{1}{k_{ex}} \int_{U_r} \frac{\xi}{\rho_\ep^2} \abs{\nab h'_\ep}^2(t) \rightarrow \int_{U_r} \xi \abs{\nab h'_*}^2(t)
\end{equation}
and hence that
\begin{equation}
\frac{1}{k_{ex}} \int_{U_r} \xi \abs{\nab h'_\ep - \nab h'_*}^2(t) \rightarrow 0.
\end{equation}
Now take $\xi$ so that $0\le \xi \le 1$ and $\xi =1$ on $U_{s_1}$ to see that
\begin{equation}
\frac{1}{k_{ex}} \int_{U_{s_1}}  \abs{\nab h'_\ep - \nab h'_*}^2(t) \rightarrow 0 \text{ for a.e. } t\in[t_1,t_2].
\end{equation}
We cannot apply the standard dominated convergence theorem directly since we lack a function that dominates the term $\int_{U_r} \abs{\dt B_\ep/\lep}^2$ in \eqref{h_con_1}.  However, 
\begin{equation}
\frac{1}{k_{ex}}\int_{U_{s_1}} \abs{\nab h'_\ep}^2 \le C + \frac{C}{k_{ex}} \int_{U_r} \abs{\frac{\dt B_\ep}{\lep}   }^2
\end{equation}
and 
\begin{equation}
\int_{t_1}^{t_2} \left(C + \frac{C}{k_{ex}} \int_{U_r} \abs{\frac{\dt B_\ep}{\lep}   }^2\right) \rightarrow \int_{t_1}^{t_2} C,
\end{equation}
so by a  variant of the dominated convergence theorem (cf. section 1.3 of \cite{ev_gar}),
\begin{equation}
\frac{1}{k_{ex}} \int_{t_1}^{t_2} \int_{U_{s_1}} \abs{\nab h'_\ep - \nab h'_*}^2 \rightarrow 0.
\end{equation}

\end{proof}

We are now in a position to present the

\begin{proof}[Proof of Theorem \ref{tensor_converge}]
 We have that 
\begin{equation}
 \nab^\bot h'_\ep + (iv_\ep,\nab_{B_\ep}v_\ep) = \frac{\dt B_\ep}{\lep} + (1-\abs{v_\ep}^2) Z_\ep \rightarrow 0
\end{equation}
in $L^2(U_r \times [t_1,t_2])$.  When combined with Lemmas \ref{tensor_swap} and \ref{h_converge}, this proves that 
\begin{equation}
 \frac{T_\ep}{k_{ex}} \rightarrow \nab^\bot h'_* \otimes \nab^\bot h'_* - I_{2\times2} \left( \hal\abs{\nab h'_*}^2 -\hal {h'_*}^2 \right)  
\end{equation}
in $L^1(U_{s_1} \times [t_1,t_2])$.  The result follows by noting that for any function $\xi:\Rn{2} \rightarrow \Rn{}$,
\begin{equation}
 \nab \xi \otimes \nab \xi + \nab^\bot \xi \otimes \nab^\bot \xi  = I_{2\times 2} \abs{\nab \xi}^2.
\end{equation}

\end{proof}

We conclude this section with a lemma that links the tensor $S$ to the renormalized energy in the case $j_{ex} =h_{ex}=1$.  The result is an adaptation of Application 3 in Chapter 13 of \cite{ss_book}.

\begin{lem}\label{h_tensor_integral}
Let $h$ solve
\begin{equation}
\begin{cases}
-\Delta h + h = 2\pi \sum_{i=1}^n d_i \delta_{a_i} & \text{in } \Omega \\
h = 0 & \text{on }\partial \Omega
\end{cases}
\end{equation}
and define the tensor 
\begin{equation}
S(h) =  -\nab h \otimes \nab h + I_{2\times2} \left(\hal \abs{\nab h}^2 + \hal \abs{h}^2 \right).
\end{equation}
Then for $r>0$ sufficiently small so that $a_j \notin B(a_i,r)$ for $j\neq i$, it holds that
\begin{equation}
\int_{\partial B(a_i,r)} S(h) \nu =  \nab_{a_i} W_{d}(a) + o_r(1),
\end{equation}
where $d =(d_1,\dotsc,d_n)$, $a=(a_1,\dotsc,a_n)$, $W_{d}(a)$ is the renormalized energy defined by \eqref{ren_def}, and $o_r(1)$ means a quantity vanishing as $r \rightarrow 0$.

\end{lem}
\begin{proof}
We begin by expanding $S(h) \nu$ in the orthonormal basis $(\nu,\tau)$ on $\partial B(a_i,r)$.  This yields
\begin{equation}
S(h)\nu = \hal \left((\nab h \cdot \tau)^2 - (\nab h \cdot \nu)^2  + h^2 \right) \nu + (\nab h \cdot \nu)(\nab h \cdot \tau) \tau.
\end{equation}
Now we use the decomposition
\begin{equation}
h =  \sum_{j=1}^n d_j(S_\Omega(\cdot,a_j) - \log\abs{\cdot -a_j}),
\end{equation}
where $S_\Omega$ is defined by \eqref{S_def}.  This allows us to write
\begin{equation}
\nab h(x) =   \sum_{j=1}^n d_j\left(\nab S_\Omega(x,a_j)  - \frac{x-a_j}{\abs{x-a_j}^2}\right) =:  -d_i \frac{x-a_i}{\abs{x-a_i}^2} + \nab H_i(x).
\end{equation}
Note that the function $H_i\in C^{1,1/2}(B(a_i,r))$.  Then on $\partial B(a_i,r)$, we have that
\begin{equation}
\nab h \cdot \nu = -\frac{d_i}{r} + \nab H_i \cdot \nu \text{ and } \nab h \cdot \tau = \nab H_i \cdot \tau,
\end{equation}
which implies that
\begin{equation}
S(h)\nu = \hal \left( -\frac{d_i^2}{r^2} + 2 \frac{\nab H_i \cdot \nu}{r}\right) \nu + \left(\frac{\nab H_i \cdot \tau}{r} \right) \tau + O(1).
\end{equation}
Integrating, we get
\begin{equation}
\int_{\partial B(a_i,r)} S(h)\nu = \frac{1}{r}\int_{\partial B(a_i,r)} \nab H_i + o_r(1)
= 2\pi \nab H_i(a_i) + o_r(1).
\end{equation}
It is then straightforward to check that $\nab H_i = \nab_{a_i} W_{d}(a),$ and the conclusion follows.
\end{proof}

%%%%%%%%%%%%%%%%%%%%%%%%%%%%%%%%%%%%%%%%%%%%%%%%%%%%%%%%%%%%%%%%%%%%%%%%
\subsection{Dynamical law}
%%%%%%%%%%%%%%%%%%%%%%%%%%%%%%%%%%%%%%%%%%%%%%%%%%%%%%%%%%%%%%%%%%%%%%%%

We have now established all of the preliminary convergence results necessary to derive the dynamical law.  As the first order of business, we calculate the divergence of the stress-energy tensor, $T_\ep$.

\begin{lem}\label{tensor_diverge}
Let $T_\ep$ be the stress-energy tensor associated to $(v_\ep,B_\ep)$, defined by \eqref{stress_def}, and let $\diverge T_\ep$ be its divergence vector, as defined by \eqref{divergence_def}.  Then
 \begin{equation}
  \diverge{T_\ep} = \frac{1}{\lep} \left( (\dt v_\ep,\nab_{B_\ep} v_\ep) - h'_\ep \dt B_\ep^\bot  \right) + (v_\ep, \nab_{B_\ep} v_\ep) \abs{Z_\ep}^2 
-  \mu(v_\ep,B_\ep) Z_\ep^\bot.
 \end{equation}
\end{lem}

\begin{proof}
A straightforward calculation shows that $\diverge{T_\ep} = (K_\ep,\nab_{B_\ep} v_\ep) -h'_\ep X_\ep^\bot$, where
\begin{equation}
K_\ep := \frac{\dt v_\ep}{\lep} -2i\nab_{B_\ep} v_\ep \cdot Z_\ep + v_\ep   \abs{Z_\ep }^2
\text{ and }
X_\ep := \frac{\dt B_\ep}{\lep} +(1-\abs{v_\ep}^2) Z_\ep. 
\end{equation}
The result follows by plugging $K_\ep,X_\ep$ into the formula for $\diverge{T_\ep}$ and noting that for any $Z\in \Rn{2}$, 
\begin{equation}
 (-2i \nab_{B_\ep} v_\ep \cdot Z,\nab_{B_\ep} v_\ep) = 2 (\partial_1^{B_\ep} v_\ep, i \partial_2^{B_\ep} v_\ep) Z^\bot
\end{equation}
and that
\begin{equation}
 (\abs{v_\ep}^2-1) h'_\ep + 2 (\partial_1^{B_\ep} v_\ep, i \partial_2^{B_\ep} v_\ep)  = - \curl(iv_\ep,\nab_{B_\ep} v_\ep) -h'_\ep = -\mu(v_\ep,B_\ep).
\end{equation}

\end{proof}

We now record a local version of the energy evolution equation  that  contains the divergence of the stress-energy tensor.

\begin{lem}\label{loc_en_evolve}
 Fix $\phi \in C_c^\infty(\Omega)$.  Then
\begin{multline}
 \dt \int_\Omega \phi \frac{\tilde{g}^0_\ep(v_\ep,B_\ep)}{\ale} + \frac{1}{\lep \ale }\int_\Omega \phi \left(\abs{\dt v_\ep}^2 + \abs{\dt B_\ep}^2 \right)   
= \frac{1}{\lep k_{ex}} \int_\Omega \phi V(v_\ep,B_\ep,0) \cdot Z_\ep \\
- \frac{1}{k_{ex}} \int_\Omega \nab \phi \cdot \left(\diverge{T_\ep} - \nab{\abs{v_\ep}^2}\frac{ \abs{Z_\ep}^2}{2}  + \mu(v_\ep,B_\ep) Z_\ep^\bot \right).
\end{multline}
\end{lem}
\begin{proof}
Note that 
\begin{equation}
\diverge(\dt v_\ep,\nab_{B_\ep} v_\ep) + \curl(h'_\ep \dt B_\ep) = \diverge( (\dt v_\ep,\nab_{B_\ep} v_\ep) - h'_\ep \dt B_\ep^\bot ). 
\end{equation}
Using this, we rewrite \eqref{accelerated_en_evolve} using Lemma \ref{tensor_diverge}, multiply by $\phi/\ale$, and integrate by parts.
\end{proof}

The following result is the main ingredient in deriving the dynamical law for the vortices.  It combines all of the previous convergence results in the case when we know that a vortex path is contained in a given cylinder and the energy is bounded in a wider cylinder.

\begin{prop}\label{local_tensor_limit}
 Suppose $a_j(t) \in B(x_0,r_1)$ for all $t\in[t_1,t_2]$ and that 
\begin{equation}\label{l_t_l_0}
 \int_{B(x_0,r_2)\backslash B(x_0,r_1)} \tilde{g}^0_\ep(v_\ep,B_\ep)(t) \le C k_{ex}
\end{equation}
for all $t\in[t_1,t_2]$.  Then for any $\phi \in C_c^\infty(\Omega)$ with the property that 
\begin{equation*}
\overline{\supp(D^2\phi)} \subset \{x\in \Omega \;\vert\; r_1 < \abs{x-x_0} < r_2 \}
\end{equation*}
it holds that 
\begin{multline}\label{l_t_l_00}
 \pi(\phi(a_j(t_2)) - \phi(a_j(t_1))) \\
= \int_{t_1}^{t_2}\int_{B(x_0,r_2)\backslash B(x_0,r_1)} D^2\phi : S(h_*) - 2\pi d_j \int_{t_1}^{t_2} Z^\bot(a_j(t)) \cdot \nab \phi(a_j(t)) dt, 
\end{multline}
where $Z:\Omega \rightarrow \Rn{2}$ is given by 
\begin{equation}
 Z = \lim_{\ep \rightarrow 0} \frac{Z_\ep}{k_{ex}} .
\end{equation}

\end{prop}
\begin{proof}
 We integrate the result of Lemma \ref{loc_en_evolve} from $t_1$ to $t_2$ to find
\begin{equation}
 I + II = III + IV + V
\end{equation}
for 
\begin{equation*}
 I:= \int_\Omega \phi \frac{\tilde{g}^0_\ep(v_\ep,B_\ep)(t_2)}{\ale} - \int_\Omega \phi \frac{\tilde{g}^0_\ep(v_\ep,B_\ep)(t_1)}{\ale},
\end{equation*}
\begin{equation*}
 II:= \frac{1}{\lep \ale}\int_{t_1}^{t_2}  \int_\Omega \phi \left(\abs{\dt v_\ep}^2 + \abs{\dt B_\ep}^2 \right),
\end{equation*}
\begin{equation*}
 III:=  \int_{t_1}^{t_2} \int_\Omega \frac{\phi}{\lep k_{ex}} V(v_\ep,B_\ep,0) \cdot Z_\ep,
\end{equation*}
\begin{equation*}
 IV:= -\frac{1}{k_{ex}} \int_{t_1}^{t_2}  \int_\Omega \nab \phi \cdot \left(\diverge{T_\ep} - \nab{\abs{v_\ep}^2}\frac{ \abs{Z_\ep}^2}{2} \right),
\end{equation*}
and
\begin{equation*}
 V:= -\int_{t_1}^{t_2}  \int_\Omega \mu(v_\ep,B_\ep) \frac{1}{k_{ex}} Z_\ep^\bot \cdot \nab \phi.
\end{equation*}
We will pass to the limit $\ep \rightarrow 0$ in each term.

From Proposition \ref{density_convergence}, we know that
\begin{equation}
 I \rightarrow \pi(\phi(a_j(t_2)) - \phi(a_j(t_1))).
\end{equation}
The bound \eqref{acc_cond2}, together with $\lep \rightarrow \infty$ implies that $ II \rightarrow 0.$  The convergence of $V(v_\ep,B_\ep,0)$ and the multiplication by $1/\lep$ show that $III \rightarrow 0$.  To handle $IV$ we first note that 
\begin{multline}
  \frac{1}{k_{ex}} \int_\Omega \nab \phi \cdot \nab{\abs{v_\ep}^2}\frac{ \abs{Z_\ep}^2}{2} = \frac{1}{2k_{ex}} \int_\Omega \nab \phi \cdot \nab{(\abs{v_\ep}^2-1)}\abs{Z_\ep}^2 \\
=\frac{1}{2k_{ex}} \int_\Omega (1-\abs{v_\ep}^2)\diverge\left( \nab \phi \abs{Z_\ep }^2\right)=o(1).
\end{multline}
By assumption, the support of $D^2\phi$ is contained in the interior of $B(x_0,r_2)\backslash B(x_0,r_1)$, and hence
\begin{equation}
 -\int_{t_1}^{t_2}  \int_\Omega \nab \phi \cdot \diverge{T_\ep} = \int_{t_1}^{t_2}  \int_{B(x_0,r_2)\backslash B(x_0,r_1)} D^2\phi : T_\ep
\end{equation}
Using the energy bound \eqref{l_t_l_0}, Theorem \ref{tensor_converge}, and a covering argument, we deduce that 
\begin{equation}
 \frac{T_\ep}{k_{ex}} \rightarrow S(h'_*) \text{ in } L^1(\overline{\supp(D^2\phi)} \times [t_1,t_2]).
\end{equation}
Hence
\begin{equation}
 IV \rightarrow \int_{t_1}^{t_2}  \int_\Omega D^2\phi : S(h'_*).
\end{equation}
For $V$ we use the convergence $\mu(v_\ep,B_\ep) \rightarrow 2\pi \sum d_i \delta_{a_i(t)}$ in conjunction 
% item 3 says that we can write $<\mu,\nab^\bot f \cdot \nab \phi>$ as the integral
with item $3$ of Proposition \ref{prod_est} to get
\begin{equation}
 V \rightarrow - 2\pi d_j \int_{t_1}^{t_2} Z^\bot(a_j(t)) \cdot \nab \phi(a_j(t)) dt.
\end{equation}

\end{proof}

We now derive the dynamical law by using this result with an appropriate test function.

\begin{thm}\label{dynamical_law}
The vortex trajectories $a_i(t),$ $i=1,\dotsc,n$, are differentiable. 
\begin{enumerate}
 \item If $h_{ex} = j_{ex} =1$, then the trajectories satisfy the dynamical law
\begin{equation}\label{d_l_0}
 \dot{a}_i(t) = -\frac{1}{\pi} \nab_{a_i}W_{d}(a(t)) -2d_i (\nab h_0(a_i(t)) - \nab^\bot f_0(a_i(t)) +  \nab^\bot f_1 (a_i(t)))
\end{equation}

 \item In the other three parameter regimes \eqref{regimes} we define 
\begin{equation}
 \lim_{\ep \rightarrow 0} \frac{j_{ex}}{k_{ex}} = \alpha \in [0,\infty) \text{ and } \lim_{\ep \rightarrow 0} \frac{h_{ex}}{k_{ex}} = \beta \in[0,\infty).
\end{equation}
Then the trajectories satisfy the dynamical law
\begin{equation}\label{d_l_00}
 \dot{a}_i(t) = -2 d_i  \alpha \nab^\bot f_1(a_i(t))  -2d_i \beta  \left(\nab h_0(a_i(t)) -\nb f_0(a_i(t))   \right).
\end{equation}
\end{enumerate}
\end{thm}

\begin{proof}
Fix a time $t_1 \in (0,T_*)$, an index $i\in\{1,\dotsc,n\}$, and a constant $0< \delta < T_*-t_1$.  By the H\"{o}lder continuity of $a_i$ and the energy bound \eqref{d_l_1}, we may find  $0<r_1<r_2 <\sigma_*/4$ with $r_1 = r_1(\delta),r_2 = r_2(\delta)$ so that $a_i(t) \in B(a_i(t_1),r_1)$ and
\begin{equation}
\int_{B(a_i(t_1),r_2) \backslash B(a_i(t_1),r_1)} \tilde{g}^0_\ep(v_\ep,B_\ep)(t) \le k_{ex} C(\delta,n,\Omega)
\end{equation}
for all $t\in[t_1,t_1+\delta]$.  Now fix a unit vector $e\in \Rn{2}$ and a function $\psi\in C_c^\infty(B(a_i(t_1),r_2))$ so that $0\le \psi \le 1$ and $\psi =1$ on $B(a_i(t_1),(r_1+r_2)/2)$.  Define the function $\phi\in C_c^\infty(\Omega)$ by $\phi(x) =  (e\cdot x)  \psi(x)$, and note that $\phi(a_i(t)) = a_i(t) \cdot e$ and $\nab \phi(a_i(t)) = e$ for all $t\in[t_1,t_1+\delta]$.  We may then apply Proposition \ref{local_tensor_limit} to deduce that for any $t_2 \in [t_1,t_1+\delta]$,
\begin{multline}
\pi (a_i(t_2) - a_i(t_1))\cdot e = - 2\pi d_i \int_{t_1}^{t_2} e\cdot Z^\bot (a_i(t))  dt
+ \int_{t_1}^{t_2}\int_{B(a_i(t_1),r_2)\backslash B(a_i(t_1),r_1)} D^2\phi : S(h'_*), 
\end{multline}
where $Z$ is defined in Proposition \ref{local_tensor_limit}.  It immediately follows that we can calculate the limit from the right:
\begin{equation}\label{d_l_3}
\pi \lim_{t \rightarrow t_1^{+}} \frac{ (a_i(t) - a_i(t_1))\cdot e}{t-t_1} = - 2\pi d_i  e\cdot Z^\bot (a_i(t_1)) + \int_{B(a_i(t_1),r_2)\backslash B(a_i(t_1),r_1)} D^2\phi : S(h'_*(t_1)).
\end{equation}
A similar argument shows that the left limit exists and agrees with the right limit. 

In the annulus  $\mathcal{A}:=B(a_i(t_1),r_2)\backslash B(a_i(t_1),r_1)$ the function $h_*(t_1)$ is smooth and satisfies $-\Delta h'_*(t_1) + h'_*(t_1) = 0$.  A direct calculation shows that for any  smooth $h$, 
\begin{equation}
\diverge S(h) = (-\Delta h+ h)\nab h, 
\end{equation}
so  $\diverge S(h'_*(t_1))=0$ in $\mathcal{A}$.  Integrating by parts and using the structure of $\phi$, this implies that 
\begin{equation}\label{d_l_4}
\int_{\mathcal{A}} D^2 \phi : S(h'_*(t_1)) = -\int_{\partial B(a_i(t_1),r_1)} e \cdot S(h'_*(t_1)) \nu.  
\end{equation}
We may thus combine \eqref{d_l_3} (with $t\rightarrow t_1^{+}$ replaced with the full limit $t \rightarrow t_1$) and \eqref{d_l_4} to deduce that $a_i$ is differentiable at $t_1$ and that 
\begin{equation}\label{d_l_5}
\pi \dot{a}_i(t_1) = - 2\pi d_i  Z (a_i(t_1))  -\int_{\partial B(a_i(t_1),r_1)}  S(h'_*(t_1)) \nu.
\end{equation}

The second through fourth parameter regimes are easier to deal with since $k_{ex} \gg 1$, which forces $h'_*=0$ and $S(h'_*)=0$.  To derive the dynamics, we must only determine the structure of $Z^\bot$ in terms of the regime.  It is easy to see that 
\begin{equation}
 Z^\bot = 
\begin{cases}
  \nb f_1 & \text{in regime }2 \\
  \nab h_0 -\nb f_0   & \text{in regime }3 \\
  \alpha \nb f_1  +\beta( \nab h_0 -\nb f_0) & \text{in regime } 4,
\end{cases}
\end{equation}
from which \eqref{d_l_00} follow immediately since $t_1$ was arbitrary.  For the regime $j_{ex}=h_{ex}=1$ we note that by sending $\delta \rightarrow 0$ we may also let $r_1 \rightarrow 0$.  Applying Lemma \ref{h_tensor_integral} and taking the limit $r_1 \rightarrow 0$ in \eqref{d_l_5} then yields
\begin{equation}
 \dot{a}_i(t_1) = -\frac{1}{\pi} \nab_{a_i}W_{d}(a(t_1)) - 2 d_i Z^\bot (a_i(t_1)), 
\end{equation}
with $Z^\bot =  \nb f_1  + \nab h_0 -\nb f_0$, from which \eqref{d_l_0} follows.

\end{proof}

%%%%%%%%%%%%%%%%%%%%%%%%%%%%%%%%%%%%%%%%%%%%%%%%%%%%%%%%%%%%%%%%%%%%%%%
\appendix
%%%%%%%%%%%%%%%%%%%%%%%%%%%%%%%%%%%%%%%%%%%%%%%%%%%%%%%%%%%%%%%%%%%%%%%

%%%%%%%%%%%%%%%%%%%%%%%%%%%%%%%%%%%%%%%%%%%%%%%%%%%%%%%%%%%%%%%%%%%%%%%%%%%%%
\section{Well-posedness and regularity}\label{well_posed_section}
%%%%%%%%%%%%%%%%%%%%%%%%%%%%%%%%%%%%%%%%%%%%%%%%%%%%%%%%%%%%%%%%%%%%%%%%%%%%%

In this section we record some results on the well-posedness and a priori estimates for \eqref{tdgl}.

\begin{prop}\label{well_posed}
The time-dependent Ginzburg-Landau equations \eqref{tdgl} are well-posed for all time, and the solutions are smooth.
\end{prop}
\begin{proof}
It is a simple matter to see that $(u,A,\Phi)$ solve \eqref{tdgl} if and only if $(v,B,\Phi)$ solve \eqref{mod_eqn_1}--\eqref{mod_eqn_2}.  Using this reformulation of the problem, the problem is amenable to standard fixed-point techniques for solving semi-linear parabolic problems.  A straightforward modification of the method employed in \cite{chen} yields well-posedness.  Smoothness follows from standard bootstrapping.
\end{proof}

We also record the following $L^\infty$ bounds on $u$ and $\nab_A u$, which follow from a simple modification of Proposition 2.8 in \cite{spirn}.

\begin{lem}\label{gradient_bound}
Suppose that the initial data satisfy $\pnormspace{u_0}{\infty}{\Omega}\le 1$ and  $\pnormspace{\nab_{A_0} u_0}{\infty}{\Omega} \le C/\ep$.  Then 
\begin{equation}
 \pnormspace{u(t)}{\infty}{\Omega} \le 1 \text{ for all } t \ge 0
\end{equation}
and
\begin{equation}
 \pnormspace{\nab_{A} u(t)}{\infty}{\Omega} \le C/\ep \text{ for all } t \ge 0.
\end{equation}
\end{lem}

%%%%%%%%%%%%%%%%%%%%%%%%%%%%%%%%%%%%%%%%%%%%%%%%%%%%%%%%%%%%%%%%%%%%%%%%
\section{Static analysis of the Ginzburg-Landau energy}\label{static_section}
%%%%%%%%%%%%%%%%%%%%%%%%%%%%%%%%%%%%%%%%%%%%%%%%%%%%%%%%%%%%%%%%%%%%%%%%

In this appendix we will record some energy estimates for the static Ginzburg-Landau energy  that are useful in the analysis of the dynamics.  In particular, when we know the limit 
\begin{equation}
\mu_\ep:=\mu(u_\ep,A_\ep) = \curl((iu_\ep,\nab_{A_\ep} u_\ep)+A_\ep) \rightarrow 2\pi \sum_{i=1}^n d_i \delta_{a_i},
\end{equation}
we will derive various estimates in terms of the vortex locations and degrees.  Most of the results are variants of well-known ones, but cannot be found in the literature in the exact form we need.

In addition to the free energy $F_\ep$ we will also use the weighted free energy $F_\ep^r$, which is given by
\begin{equation}
F_\ep^r(u,A) = \hal \int_\Omega \abs{\nab_A u}^2 + \frac{1}{2\ep^2}(1-\abs{u}^2)^2 + r^2\abs{\curl{A}}^2
\end{equation}
for some $r>0$.  The simplified energy is given by $E_\ep(u) = F_\ep(u,0)$.  We employ the notation $F_\ep(u,A,S)$ for $S\subset \Omega$ to mean the energy with the integral evaluated only over $S$.

Our first result gives a lower bound for the free energy in terms of the limiting vorticity measure.

\begin{prop}\label{lower_bound}
 Suppose $(u_\ep,A_\ep)$ satisfy the bound $F_\ep(u_\ep,A_\ep) \le C \ale$ as well as $\pnorm{u_\ep}{\infty}\le 1,$ $\pnorm{\nab_{A_\ep} u_\ep}{\infty} \le C/\ep,$ and $\nab_{A_\ep} u_\ep \cdot \nu=0$ on $\partial \Omega$.
Further suppose that 
\begin{equation}
 \mu_\ep \rightarrow 2\pi \sum_{i=1}^n d_i \delta_{a_i},
\end{equation}
where $d_i = \pm 1$, and that $\pnorm{\curl{A_\ep}}{2}\le K_\ep$ for $0 \le K_\ep \ll \ale$.  Then
\begin{equation}
 F_\ep(u_\ep,A_\ep) \ge  \pi n \ale + n \gamma + W_{d}(a) + o(1)K_\ep.
\end{equation}

\end{prop}
\begin{proof}
For convenience we will drop the subscript $\ep$ in the proof, writing $u,A,\mu$ in place of $u_\ep,A_\ep,\mu_\ep$.  Since we want a lower bound of $F_\ep$, we may first perform a minimization of $F_\ep(u,A)$ over all $A,$ keeping $u$ fixed.  That is, we bound
\begin{equation}
 F_\ep(u,A) \ge \min_{A} F_\ep(u,A) := F_\ep(u,B)
\end{equation}
where $B$ solves 
\begin{equation}\label{l_b_1}
 \begin{cases}
  \nab^\bot \curl{B} = -(iu,\nab_B u) &\text{in } \Omega \\
  \curl{B} = 0 & \text{on }\partial \Omega.
 \end{cases}
\end{equation}
We also fix the Coulomb gauge so that $\diverge{B} = 0$ in $\Omega$ and $B\cdot \nu =0$ on $\partial \Omega$.
By the Poincar\'{e} lemma, this allows us to write $B = \nab^\bot \xi,$ where $\xi = 0$ on $\partial \Omega$.  Since $\curl{B} = \Delta \xi$ and $\abs{\nab \curl{B}}\le \abs{\nab_B v}$, elliptic regularity gives that 
\begin{equation}\label{l_b_3}
\begin{split}
 & \norm{\xi}_{H^2(\Omega)} \le \pnorm{\curl{A}}{2} \le K_\ep, \\
 &  \norm{\xi}^2_{H^3(\Omega)}\le C \int_\Omega \abs{\curl{A}}^2 + \abs{\nab_A v}^2 \le C F_\ep(v,A) \le C \ale.
\end{split}
\end{equation}
We see from the first bound that $\norm{\xi}_{C_0^{0,\alpha}(\Omega)} \le C K_\ep$ for any $\alpha \in (0,1)$, and from the second  that $\pnorm{A}{\infty}^2\le C \ale$.

We expand $F_\ep(u,B),$ employing the fact that $B$ solves \eqref{l_b_1}, to see that 
\begin{equation}
 F_\ep(u,B) = E_\ep(u) - \hal \int_\Omega \abs{B}^2 + \abs{\curl{B}}^2 + o(1).
\end{equation}
Taking the curl of \eqref{l_b_1} and adding $\curl{B}$ to both sides, we have that $-\Delta^2 \xi + \Delta \xi = \mu$ in $\Omega$, with $\xi=\Delta\xi=0$ on $\partial \Omega$. Then, integrating by parts, we have that 
\begin{equation}\label{l_b_4}
\begin{split} 
\hal \int_\Omega \abs{B}^2 + \abs{\curl{B}}^2 & = \hal \int_\Omega \abs{\nab \xi}^2 + \abs{\Delta \xi}^2 \\
& = \hal \int_\Omega \xi (\Delta^2 \xi- \Delta \xi)  \\
& = -\hal\int_\Omega \xi \mu.
\end{split}
\end{equation}
Define $\xi_*$ as the solution to 
\begin{equation}
\begin{cases}
-\Delta^2 \xi_* + \Delta \xi_* =   2\pi \sum d_i \delta_{a_i} & \text{in }\Omega \\
\xi_*=\Delta \xi_* =0  & \text{on }  \partial \Omega.
\end{cases}
\end{equation}
Since $\xi/K_\ep$ is bounded in $C_0^{0,\alpha}(\Omega)$ for any $\alpha\in(0,1)$, we may assume that up to extraction 
\begin{equation}
 \frac{\xi}{K_\ep} \rightarrow 
\begin{cases}
 \xi_* & \text{if } K_\ep = O(1) \\
 0 & \text{if } 1 \ll K_\ep \ll \ale.
\end{cases}
\end{equation}
In either case
\begin{equation}\label{l_b_6}
 \int_\Omega \xi \mu = 2\pi \sum_{i=1}^n d_i \xi_*(a_i) + o(1)K_\ep.
\end{equation}

To deal with the $E_\ep(u)$ term we have to show that the hypotheses carry over to $u$.  First note that 
\begin{equation}
 \pnorm{\nab u}{\infty} \le \pnorm{\nab_B u}{\infty} + \pnorm{u}{\infty} \pnorm{B}{\infty} \le \frac{C}{\ep} + C \sqrt{\ale} \le \frac{C}{\ep}
\end{equation}
and that  $E_\ep(v) \le C \ale.$
Since $B\cdot \nu=0$ on the boundary, it holds that $\nab u \cdot \nu = \nab_B v \cdot \nu =0$ on $\partial \Omega$.  Finally, we note that since $\mu = \curl( (iu,\nab u) +(1-\abs{u}^2)B),$ the convergence of $\mu$ guarantees that $ \curl(iu,\nab u) \rightarrow 2\pi \sum_{i=1}^n d_i \delta_{a_i}$ as well.  We thus have that all of the hypotheses of Proposition 4.3 of \cite{serf_1} are satisfied; we find that 
\begin{equation}\label{l_b_5}
 E_\ep(u) \ge \pi n \ale + n \gamma +\left(-\pi \sum_{i\neq j}d_i d_j \log\abs{a_i-a_j} +\pi \sum_{i,j}d_i d_j R_\Omega(a_i,a_j)    \right) +o(1),
\end{equation}
where $R_\Omega$ is defined by $ R_\Omega(x,y) = P(x,y) + \log\abs{x-y},$ with $P$ is the solution to
\begin{equation}
 \begin{cases}
  -\Delta_x P(x,y) = 2\pi \delta_{y} & \text{in } \Omega \\
  P(x,y) =0 & \text{for } x\in\partial \Omega.
 \end{cases}
\end{equation}

We deduce that 
\begin{multline}\label{l_b_2}
 F_\ep(u,A) \ge \pi n \ale + n \gamma    \\
+ \left(-\pi \sum_{i\neq j} d_i d_j \log\abs{a_i-a_j} + \pi \sum_{i=1}^n d_i \xi_*(a_i) + \pi \sum_{i,j}d_i d_j  R_\Omega(a_i,a_j)    \right) + o(1)K_\ep.
\end{multline}
A simple calculation shows that 
\begin{equation}
  \sum_{i=1}^n d_i \xi_*(a_i) + \sum_{i,j}d_i d_j R_\Omega(x,a_j) =  \sum_{i,j} d_i d_j S_\Omega(a_i,a_j).
\end{equation}
Substituting this into \eqref{l_b_2} yields the desired lower bound. 

\end{proof}

The next result gives a lower bound for the weighted free energy in balls near the vortex locations.

\begin{lem}\label{inverse_measure}
Suppose that $(u_\ep,A_\ep)$ are such that 
\begin{equation}\label{i_m_2}
 F_\ep(u_\ep, A_\ep) \le \pi n \ale + o(\ale)
\end{equation}
and that 
\begin{equation}\label{i_m_1}
\mu_\ep \rightarrow \mu = 2\pi \sum_{i=1}^n d_i \delta_{a_i},
\end{equation}
where $d_i = \pm 1$.  Let $0< \sigma < \hal \min\{ \dist(a_i,\partial \Omega)\} \cup \{\abs{a_i - a_j} \;\vert\; i\neq j  \}$.  Then there is a universal constant $C>0$ so that
\begin{equation}\label{i_m_00}
  \pi n \left(\log \frac{\sigma}{\ep} - C \right) \le F^{\sigma/3}_\ep(u_\ep,A_\ep,\cup B(a_i,\sigma)).
\end{equation}

\end{lem}

\begin{proof}
We apply the ball construction (Theorem 4.1 of  \cite{ss_book}) in each ball $B(a_i,\sigma)$ to find a collection of balls $\mathcal{B}_i = \{B_{i,1},\dotsc,B_{i,m_i}\}$ with final radius $r = \ale^{-2}$. For each $i=1,\dotsc,n$ and $j=1,\dotsc,m_i$, we let $b_{i,j}$ denote the center of the ball $B_{i,j}$ and $d_{B_{i,j}}$ denote its degree.  Then, according to the Theorem 6.1 of \cite{ss_book}, for each $i=1,\dotsc,n$ it holds that 
\begin{equation}
 \norm{\mu_\ep - 2\pi \sum_{\substack{1\le j\le m_i \\ B_{i,j} \subset B(a_i,\sigma-\ep) }} d_{B_{i,j}} \delta_{b_{i,j}} }_{(C_0^{0,1}(B(a_i,\sigma)))^{*}} \le \frac{C}{\ale^2} F_\ep(u_\ep,A_\ep) = o(1).
\end{equation}
From this and assumption \eqref{i_m_1}, we deduce that for each $i=1,\dotsc,n$ there exists a ball $B_{i,j} \in \mathcal{B}_i$ such that $B_{i,j} \subset B(a_i,\sigma-\ep)$, $\abs{d_{B_{i,j}}} \ge 1,$ and $b_{i,j} \rightarrow a_i$ as $\ep \rightarrow 0$.  By relabeling the collection $\mathcal{B}_i$ we may assume that $B_{i,j} = B_{i,1}$.

The energy estimates of the ball construction imply that 
\begin{equation}\label{i_m_3}
 \pi D_i \left(\log \frac{r}{\ep D_i} - C \right) \le F^{r}_\ep(u_\ep,A_\ep,\cup_j B_{i,j}),
\end{equation}
where
\begin{equation*}
D_i = \sum_{\substack{1\le j\le m_i \\ B_{i,j} \subset B(a_i,\sigma-\ep) }} \abs{d_{B_{i,j}}}.
\end{equation*}
Moreover, $D_i \le C n$ for some universal constant $C$.  Summing \eqref{i_m_3} over $i=1,\dotsc,n$ and plugging in the value of $r$ and the upper bound \eqref{i_m_2}, we find that
\begin{equation}
 \pi \sum_{i=1}^n D_i \ale - 2\pi \sum_{i=1}^n D_i \log \ale - \pi \sum_{i=1}^n D_i (\log D_i +C ) \le \pi n \ale + o(\ale).
\end{equation}
Dividing by $\ale$, we deduce that $\sum D_i \le n$ for $\ep$ is sufficiently small.  However, since for each $i=1,\dotsc,n$ the ball $B_{i,1} \subset B(a_i,\sigma - \ep)$ and $\abs{d_{B_{i,1}}} \ge 1$, the reverse inequality $n\le \sum D_i$ must also hold.  Hence $\sum D_i=n$.  From this we deduce that $D_i = 1$ for each $i$ and that $d_{B_{i,j}}=0$ for $j=2,\dotsc,m_i$.

We now grow the balls in each collection $\mathcal{B}_i$ into a larger collection, $\mathcal{B}_i' $, with total radius $r=\sigma/3$.  There must exist a ball $B_i' \in \mathcal{B}_i'$ such that $B_{i,1} \subset B_i'$.  By the above analysis on the degrees of the balls in $\mathcal{B}_i$ it holds that $\deg(B_i')= d_i$ and that the degrees of the other balls in $\mathcal{B}_i'$ vanish.  Moreover, since $\text{rad}(B_i') \le \sigma/3$ and $B_i' \ni b_{i,1} \rightarrow a_i$, it must be the case that $B_i' \subset B(a_i,\sigma-\ep)$.  Plugging into the ball construction again then yields the bound
\begin{equation}
  \pi \left(\log \frac{\sigma}{3\ep } - C \right) \le F^{\sigma/3}_\ep(u_\ep,A_\ep,B(a_i,\sigma))
\end{equation}
for each $i=1,\dotsc,n$, from which \eqref{i_m_00} immediately follows.
\end{proof}

We now use the previous lemmas to prove that $\curl{A_\ep}$ is lower order than $F_\ep$ and that the energy away from the vortex locations is bounded.

\begin{lem}\label{outside_bound}
 Suppose $(u_\ep,A_\ep)$ satisfy $\mu_\ep \rightarrow 2\pi \sum_{i=1}^n d_i \delta_{a_i}$ with $d_i = \pm 1$ as well as the bound
\begin{equation}
F_\ep(u_\ep,A_\ep) \le  \pi n \ale  + K_\ep
\end{equation}
for some sequence $0\le K_\ep \ll \ale$.  Let 
\begin{equation}
0 < \sigma < \frac{1}{2} \min\{ \dist(a_i,\partial \Omega)\} \cup \{\abs{a_i - a_j} \;\vert\; i\neq j  \}\cup\{6/\sqrt{2}\}.
\end{equation}
Write $\Omega_\sigma = \Omega \backslash \cup B(a_i,\sigma)$. Then 
\begin{equation}
 \hal \int_{\Omega_\sigma}\abs{\nab_{A_\ep} u_\ep}^2  + \frac{1}{2\ep^2} (1-\abs{u_\ep}^2)^2  + \frac{1}{8} \int_\Omega \abs{\curl{A_\ep}}^2 
\le  \pi n \log \frac{1}{\sigma} + Cn + K_\ep ,
\end{equation}
where $C_\Omega$ is a constant depending on $\Omega$ and $C$ is a universal constant.
\end{lem}
\begin{proof}
Write $ F_\ep(u_\ep,A_\ep) = F_\ep(u_\ep,A_\ep,\Omega_\sigma) + F_\ep(u_\ep,A_\ep,\cup B(a_i,\sigma)).$  We further decompose the latter term to
\begin{equation}
 F_\ep(u_\ep,A_\ep,\cup B(a_i,\sigma)) = F^{\sigma/3}_\ep(u_\ep,A_\ep,\cup B(a_i,\sigma)) + \hal\left(1 - \frac{\sigma^2}{9} \right) \int_{\cup B(a_i,\sigma)} \abs{\curl{A_\ep}}^2.
\end{equation}
By construction
\begin{equation}
 \hal\left(1 - \frac{\sigma^2}{9} \right) \int_{\cup B(a_i,\sigma)} \abs{\curl{A_\ep}}^2 \ge \frac{1}{4} \int_{\cup B(a_i,\sigma)} \abs{\curl{A_\ep}}^2.
\end{equation}
Lemma \ref{inverse_measure} shows that $ F^{\sigma/3}_\ep(u_\ep,A_\ep,\cup B(a_i,\sigma)) \ge \pi n \left(\log \frac{\sigma}{\ep} - C \right) $
for $C$ a universal constant.  Chaining together these lower bounds with the upper bound of the hypothesis then gives
\begin{equation}
 \hal \int_{\Omega_\sigma}\abs{\nab_{A_\ep} u_\ep}^2  + \frac{1}{2\ep^2} (1-\abs{u_\ep}^2)^2 +\frac{1}{4} \int_\Omega \abs{\curl{A_\ep}}^2 
 \le  \pi n \log \frac{1}{\sigma} +Cn + K_\ep,
\end{equation}
the desired inequality.
\end{proof}

%%%%%%%%%%%%%%%%%%%%%%%%%%%%%%%%%%%%%%%%%%%%%%%%%%%%%%%%%%%%%%%%%%%%%%%%
%References
%%%%%%%%%%%%%%%%%%%%%%%%%%%%%%%%%%%%%%%%%%%%%%%%%%%%%%%%%%%%%%%%%%%%%%%%

\end{document}